\documentclass[a4paper,11 pt]{amsart}

\RequirePackage{amsmath, amssymb, amsthm, amsfonts}
\usepackage{graphicx}
\usepackage{multicol}
\usepackage{pifont}
\usepackage{mathrsfs}
\usepackage{units}
\usepackage{xcolor}
\RequirePackage[all]{xy}
\usepackage{fullpage}
\usepackage[mathscr]{eucal}
\usepackage{float}
\newtheorem{theo}{Theorem}[section]
\newtheorem{lem}[theo]{Lemma}
\newtheorem{prop}[theo]{Proposition}
\newtheorem{cor}[theo]{Corollary}
\newtheorem{defin}[theo]{Definition}

\theoremstyle{definition}

\newtheorem{rem}[theo]{Remark}

\newcommand{\SL}{{\mathrm{SL}}}

\newcommand{\GL}{{\mathrm{GL}}}
\newcommand{\Aut}{{\mathrm{Aut}}}
\newcommand{\Out}{{\mathrm{Out}}}
\newcommand{\Inn}{{\mathrm{Inn}}}

\newcommand{\Diff}{{\mathrm{Diff}}}
\newcommand{\Map}{{\mathrm{Map}}}

\newcommand{\lr}{\longrightarrow}

\title[Isotrivial surfaces with $p_g=q=2$]{The classification of isotrivially fibred surfaces with $p_g=q=2$}
\author{ Matteo Penegini}%, Soenke Rollenske
\address{Matteo Penegini\\Lehrstuhl Mathematik VIII\\
Universit\"at Bayreuth, NWII\\
D-95440 Bayreuth, Germany} \email{matteo.penegini@uni-bayreuth.de}
\address{S\"onke Rollenske\\Lehrstuhl Mathematik VIII\\
Universit\"at Bayreuth, NWII\\
D-95440 Bayreuth, Germany}
\email{soenke.rollenske@uni-bayreuth.de}
\begin{document}
%%%%%%%%%%%%%%%%%%%%%%%%%%%%%%%%%%%%%%%%%%%%%%%%%%%%%%%%%%%%%%%%%%%%%%%%%%%%%%%%%%%%%%%%%%%
%%%%%%%%%%%%%%%%%%%%%%%%%%%%%%%%%%%%%%%%%%%%%%%%%%%%%%%%%%%%%%%%%%%%%%%%%%%%%%%%%%%%%%%%%%%
\maketitle
%%%%%%%%%%%%%%%%%%%%%%%%%%%%%%%%%%%%%%%%%%%%%%%%%%%%%%%%%%%%%%%%%%%%%%%%%%%%%%%%%%%%%%%%%%%
\begin{abstract}
An isotrivially fibred surface is a smooth projective surface
endowed with a morphism onto a curve such that all the smooth
fibres are isomorphic to each other. The first goal of this paper
is to classify the isotrivially fibred surfaces with $p_g=q=2$
completing and extending a result of Zucconi. As an important
byproduct, we provide new examples of minimal surfaces of general
type with $p_g=q=2$ and $K^2=4,5$ and a first example with
$K^2=6$.
\end{abstract}
%%%%%%%%%%%%%%%%%%%%%%%%%%%%%%%%%%%%%%%%%%%%%%%%%%%%%%%%%%%%%%%%%%%%%55%%%%%%%%%%%%%%%%%%%
%%%%%%%%%%%%%%%%%%%%%%%%%%%%%%%%%%%%%%%%%%%%%%%%%%%%%%%%%%%%%%%%%%%%%%%%%%%%%%%%%%%%%%%%%%
\begin{center}
%Versione del 01/02/2010
\end{center}
\section{Introduction}
The classification of smooth connected minimal complex projective
surfaces of general type with small invariants is far from being
achieved, and up to now a complete classification seems out of reach. This is a reason why
one first tries to understand and classify surfaces with
particularly small invariants, for example with
$\chi(\mathcal{O}_S)=1$. If this is the case one has
$1=\chi(\mathcal{O}_S)=1-q+p_g$ and it follows that $p_g=q$. If we
also assume that the surface is irregular (i.e., $q
>0$) then the Bogomolov-Miyaoka-Yau and Debarre inequalities,
$K^2_S \leq 9$, $K^2_S \geq 2p_g$, imply $1 \leq p_g \leq 4$. If $p_g=q=4$ we
have a product of curves of genus 2, as shown by Beauville, while
the case $p_g=q=3$ was understood through the work of several
authors \cite{CCML, hpar, pir}. It seems that, the
classification becomes more complicated as the value of $p_g$
decreases. In this paper we address the case  $p_g=q=2$.

We say that a surface $S$ is \emph{isogenous to a product of curves} if $S = (C \times F )/G$,
for $C$ and $F$ smooth curves and $G$ a finite group acting freely on $C \times F$. Surfaces isogenous to a product were introduced by Catanese in \cite{cat00}. They are of general type if and only if both $g(C)$
and $g(F)$ are greater than or equal to $2$ and in this case $S$ admits a unique minimal realization where they are as small as possible. From now on, we tacitly assume that such a realization is
chosen, so that the genera of the curves and the group $G$ are invariants of $S$.
%The action of $G$
%can be seen to respect the product structure on $C \times F$ . This means that such actions fall in
We have two cases: the \emph{mixed} one, where there exists some
element in $G$ exchanging the two factors (in this situation $C$
and $F$ must be isomorphic) and the \emph{unmixed} one, where $G$
acts faithfully on both $C$ and $F$ and diagonally on their
product. A special case of surfaces isogenous to a product of
unmixed type is the case of \emph{generalized hyperelliptic
surfaces} where $G$ acts freely on $C$ and $F/G \cong \mathbb{P}^1$.

A generalization of the unmixed case is the following: consider a
finite group $G$ acting faithfully on two smooth projective curves
$C$ and $F$ of genus $\geq 2$, and diagonally, but not necessarily
freely, on their product, and take the minimal resolution
$S'\rightarrow T:=(C \times F)/G$ of the singularities of $T$. In
this case the holomorphic map:
\[ f'_1\colon S' \longrightarrow C':=C/G
\]
is called a \emph{standard isotrivial fibration} if it is a
relatively minimal fibration. More generally an \emph{isotrivial fibration} is a fibration $f:S \rightarrow B$  
from a smooth surface onto a smooth curve such that all the smooth fibres are isomorphic to each other. A monodromy argument shows that, 
in case the general fibre $F$ is irrational, there is a birational realization of $S$ as a quotient of a product of two curves $S \stackrel{bir}{\sim}(C \times F)/G \rightarrow C/G\cong B$.

Among the surfaces which admit an isotrivial fibration one can find examples of surfaces with
$\chi(\mathcal{O}_S)=1$. Since \cite{cat00} appeared several
authors studied intensively standard isotrivially fibred
surfaces, and eventually classified all those, which are minimal with $p_g=q=0$ \cite{infabfriz, ingefabfritzroby, BP} and with $p_g=q=1$ \cite{pol2, canpol, mistretun}.

In this paper we complete the classification of isotrivially fibred
surfaces with $p_g=q=2$ (which was partially given in
\cite{zucc}).
Moreover we give a precise description of the
corresponding locus in the moduli space of surfaces of general
type. Indeed, by the results of \cite{cat00}, this locus is a
union of connected components in the case of surfaces isogenous to a product of curves,
and irreducible subvarieties in the case of only isotrivial surfaces.
We calculate the number of these
components (subvarieties) and their dimensions. The following theorem summarizes
our classification.
%---------------------------------------------------------------------------------------------------------------%
\begin{theo}\label{zerozero}
Let $S$ be a minimal surface of general type with
$p_g=q=2$ such that it is either a surface isogenous to a product of curves of mixed type or it admits an isotrivial fibration. Let $\alpha \colon S \rightarrow Alb(S)$ be the Albanese map.
Then we have the following possibilities:
\begin{enumerate}
\item If $dim(\alpha(S))=1$, then $S \cong (C \times F)/G$ and it
is generalized hyperelliptic. The classification of these
surfaces is given by the cases labelled with GH in Table 1,
where we specify the only possibilities for the genera of the
two curves $C$ and $F$, and for the group $G$.

\item If $dim(\alpha(S))=2$, then there are three cases:
\begin{itemize}
\item $S$ is isogenous to product of unmixed type $(C \times
F)/G$, and the classification of these surfaces is given by the
cases labelled with UnMix in Table 1;

\item $S$ is isogenous to a product of mixed type $(C \times
C)/G$, there is only one case and it is given in Table 1
labelled with Mix;

\item $S \rightarrow T:=(C \times F)/G$ is a minimal
desingularization of $T$, and these surfaces are classified in
Table 2.
\end{itemize}
\end{enumerate}
\begin{center}
\begin{tabular}{|c|c|c|c|c|c|c|c|c|}
  \hline
 Type & $K^2_S$ & $g(F)$ & $g(C)$ & $G$ & IdSmallGroup & \textbf{m} & $dim$ & $n$ \\
  \hline
 GH & $8$ & $2$ & $3$ &  $\mathbb{Z}_2$ &                    G(2,1) & $(2^6)$ & $6$ & $1$ \\ \hline
 GH & $8$ &   $2$ & $4$ &  $\mathbb{Z}_3$ &                     G(3,1) & $(3^4)$ & $4$ & $1$ \\ \hline
 GH & $8$ &   $2$ & $5$ &  $\mathbb{Z}_2 \times \mathbb{Z}_2$ & G(4,2) & $(2^5)$ & $5$ & $2$  \\ \hline
 GH & $8$ &   $2$ & $5$ &  $\mathbb{Z}_4$ &                     G(4,1) & $(2^2,4^2)$ & $4$ & $1$ \\ \hline
 GH & $8$ &   $2$ & $6$ &  $\mathbb{Z}_5$ &                     G(5,1) & $(5^3)$ & $3$ & $1$ \\ \hline
 GH & $8$ &   $2$ & $7$ &  $\mathbb{Z}_6$ &                     G(6,2)  & $(2^2,3^2)$ & $4$ & $1$\\ \hline
 GH & $8$ &  $2$ & $7$ &  $\mathbb{Z}_6$ &                     G(6,2) & $(3,6^2)$ & $3$ & $1$\\ \hline
 GH & $8$ &   $2$ & $9$ &  $\mathbb{Z}_8$ &                     G(8,1) & $(2,8^2)$ & $3$ & $1$\\ \hline
%  \end{tabular}
% \end{center}
%\begin{center}
%\begin{tabular}{|c|c|c|c|c|c|c|c|c|}
%  \hline Type & $K^2_S$ & $g(F)$ & $g(C)$ & $G$ & IdSmallGroup & \textbf{m} & $dim$ & $n$ \\
%  \hline
 GH & $8$ &  $2$ & $11$ & $\mathbb{Z}_{10}$&                   G(10,2) & $(2,5,10)$ & $3$ & $1$ \\ \hline
 GH & $8$ &  $2$ & $13$ & $\mathbb{Z}_2 \times \mathbb{Z}_6$ & G(12,5) & $(2,6^2)$ & $3$ & $2$\\ \hline
 GH & $8$ &  $2$ & $7$ &  $S_3$ &                              G(6,1) & $(2^2,3^2)$ & $4$ & $1$\\ \hline
 GH & $8$ &  $2$ & $9$ &  $Q_8$ &                              G(8,4) & $(4^3)$ & $3$ & $1$ \\ \hline
 GH & $8$ &  $2$ & $9$ &  $D_4$ &                              G(8,3) & $(2^3,4)$ & $4$ & $2$ \\ \hline
 GH & $8$ &  $2$ & $13$ & $D_6$ &                              G(12,4) & $(2^3,3)$ & $3$ & $2$ \\ \hline %WW
 GH & $8$ &  $2$ & $13$ & $D_{4,3,-1}$ &                       G(12,1) & $(3,4^2)$ & $3$ & $1$ \\ \hline
 GH & $8$ &  $2$ & $17$ & $D_{2,8,3}$ &                        G(16,8) & $(2,4,8)$ & $3$ & $1$ \\ \hline
 GH & $8$ &  $2$ & $25$ & $\mathbb{Z}_2 \ltimes ((\mathbb{Z}_2)^2 \times \mathbb{Z}_3)$ & G(24,8) & $(2,4,6)$ & $3$ & $2$ \\ \hline
 GH & $8$ &  $2$ & $25$ & $\SL(2,\mathbb{F}_3)$ &               G(24,3) & $(3^2,4)$ & $3$ & $1$ \\ \hline
 GH & $8$ &  $2$ & $49$ & $GL(2,\mathbb{F}_3)$ &               G(48,29) & $(2,3,8)$ & $3$ & $1$ \\ \hline
UnMix & $8$ & $3$ & $3$ & $\mathbb{Z}_2 \times \mathbb{Z}_2$ & G(4,2) & $(2^2)$, $(2^2)$ & $4$ &$1$  \\ \hline
UnMix & $8$ &  $3$ & $4$ & $S_3$                              & G(6,1) & $(3)$, $(2^2)$ & $3$ &$1$  \\ \hline
UnMix & $8$ &  $3$ & $5$ & $D_4$ &                              G(8,3) & $(2)$, $(2^2)$ & $3$ & $1$ \\ \hline
Mix & $8$ & $3$ & $3$  & $\mathbb{Z}_4$ &                      G(4,1) & - & $3$ & $1$ \\ \hline
\end{tabular}
\end{center}
\begin{center}
Table 1.
\end{center}
In Table 1 and 2 \texttt{IdSmallGroup} denotes the label of the group $G$ in the GAP4
database of small groups, $\mathbf{m}$ is the branching data. In Table 1 each item provides a union of connected components of the moduli space of surfaces of general
type, their dimension is listed in the column $dim$ and $n$ is the
number of connected components.
\begin{center}
\begin{tabular}{|c|c|c|c|c|c|c|c|c|c|}
  \hline
 $K^2_S$ & $g(C)$ & $g(F)$ & $G$ & \texttt{IdSmallGroup} & \textbf{m} & Type & Num.~Sing.  & dim & n
\\ \hline
$4$ &  $2$ & $2$ &  $\mathbb{Z}_2$  & G(2,1)  & $(2^2)$ $(2^2)$&
$\frac{1}{2}(1,1)$ & $4$ & $4$ & $1$
\\ \hline $4$ &  $3$ & $3$ &  $D_4$ & G(8,3)  & $(2)$ $(2)$ & $\frac{1}{2}(1,1)$ & $4$ & $2$ & $1$
\\ \hline $4$ & $3$ & $3$ &  $Q_8$ & G(8,4)  & $(2)$ $(2)$ & $\frac{1}{2}(1,1)$ & $4$ & $2$ & $1$
\\ \hline $5$ &  $3$ & $3$ &  $S_3$ & G(6,1) &  $(3)$ $(3)$ & $\frac{1}{3}(1,1)+\frac{1}{3}(1,2)$ & $2$ & $2$ & $1$
\\ \hline $6$ &  $4$ & $4$ &  $A_4$ & G(12,3) & $(2)$ $(2)$ & $\frac{1}{2}(1,1)$ & $2$ & $2$ & $1$
\\
\hline
\end{tabular}
Table 2.
\end{center}
In Table 2 each item provides a union of irreducible subvarieties of the moduli space of surfaces of general type, their dimension is listed in the column $dim$ and $n$ is the
number of subvarieties. Moreover the columns of Table 2 labelled with $Type$ and $Num.~Sing.$ indicate the types and the number of
singularities of $T$.
\end{theo}
%---------------------------------------------------------------------------------------------------------------------%

We point out that in Table 2 there are new examples of minimal
surfaces of general type with $p_g=q=2$ and $K^2_S=4,5$, and a
first example with $K^2_S=6$. It will be interesting to find, if
there are any, examples of surfaces with $p_g=q=2$ and
$K^2_S=7$ or $9$.

We recall that surfaces of general type with
$p_g=q=2$ and $K^2_S=4$ were studied by Ciliberto and Mendes
Lopes. Indeed they proved that the surfaces with $p_g=q=2$ and non-birational
bicanonical map are double coverings of a principally polarized
abelian surfaces branched on a divisor $D \in |2 \Theta|$, and
they have $K^2_S=4$ (\cite{cilimen}). While Chen and Hacon
(\cite{chenha}) constructed a first example of a surface with
$K^2_S=5$.

The classification of the surfaces isogenous to a
product involves techniques from both geometry and combinatorial group
 theory (developed
in \cite{ingridfab, infabfriz, pol1, canpol}). In the second section of this paper we recall the
 relation between coverings of Riemann surfaces and orbifold surface groups,
 which enables us to transform the geometric problem of classification into an
 algebraic one.

 In the third section we recall the notion of
 generalized hyperelliptic surface, and we shall see, following
 \cite{cat00} and \cite{zucc}, that all the surfaces with $p_g=q=2$ and not of Albanese general type
are generalized hyperelliptic.
 Using this fact and the material of section two we
 classify all such surfaces. We notice that such
 classification was partially given in \cite{zucc} using different techniques.
We proceed then to classify the surfaces isogenous to a product of
unmixed type and of Albanese general type and finally we
study the mixed case.

In the fourth section we consider the case when $(C \times
F)/G$ is singular. We provide there new examples of surfaces with
$p_g=q=2$. 

 In the fifth section we treat the problem of the description of the
 moduli spaces. We recall a theorem of Bauer and Catanese
 (Theorem 1.3 \cite{ingridfab} ) which
 tells how to calculate the number of the connected components of the moduli
 space. Moreover we calculate the action of the mapping class
 group of a Riemann surface of genus $2$ on its fundamental
 group. In the end we shall see that calculating the number of
 connected components is a task which cannot be achieved easily
 without using a computer.

In the sixth section we calculate the fundamental groups of the found isotrivial surfaces.
We recall in that section two structure theorems one
for the fundamental group of surfaces isogenous to a product of curves and one for the fundamental group of isotrivial surfaces
following \cite{ingefabfritzroby}.

 In the appendix written by S. Rollenske, it is described the GAP4 program we wrote
 to finish our classification and it is explained how to use it. Since
 it is written in great generality we hope it can be used for
 other tasks. For example we were able to finish the classification given
 in \cite{canpol} adding the number of connected components of
 the moduli spaces. We tested the program also on the cases with
 $p_g=q=0$ given in \cite{infabfriz} and $p_q=q=1$ given in \cite{pol1}.
%%%%%%%%%%%%%%%%%%%%%%%%%%%%%%%%%%%%%%%%%%%%%%%%%%%%%%%%%%%%%%%%%%%%%%%%%%%%%%%%%%%%%%%%%%%%%%%%%%%%%%%%%%%%%%%%%%%%%%%%%%%%
%%%%%%%%%%%%%%%%%%%%%%%%%%%%%%%%%%%%%%%%%%%%%%%%%%%%%%%%%%%%%%%%%%%%%%%%%%%%%%%%%%%%%%%%%%%%%%%%%%%%%%%%%%%%%%%%%%%%%%%%%%%%

\bigskip

\textbf{Acknowledgements.} The author would like to thank especially Fabrizio Catanese for suggesting this research,
Francesco Polizzi for pointing out the incompleteness of the
classification of isotrivial surfaces with $p_g=q=2$ and many
suggestions, and S\"onke Rollenske for his contribution. Moreover
he thanks I. Bauer, M. L\"onne, E. Mistretta, R. Pignatelli for
useful discussion and suggestions.

The author acknowledges the support of the DFG Forschergruppe 790
''Classification of algebraic surfaces and compact complex
manifolds''.
%%%%%%%%%%%%%%%%%%%%%%%%%%%%%%%%%%%%%%%%%%%%%%%%%%%%%%%%%%%%%%%%%%%%%%%%%%%%%%%%%%%%%%%%%%%%%%%%%%%%%%%%%%%%%%%%%%%%%%%%%%%%%%

\bigskip 

\textbf{Notation and conventions.} 
We shall denote by $S$ a smooth, irreducible,
complex, projective surface. We shall also use the
standard notation in surface theory, hence we denote by $\Omega^p_S$ the sheaf of holomorphic $p-$forms on $S$,
$p_g:=h^0(S,\Omega^2_S)$ the \emph{geometric genus} of $S$,
$q:=h^0(S,\Omega^1_S)$ the \emph{irregularity} of $S$,
$\chi(S)=\chi(\mathcal{O}_S)=1+p_g-q$ the \emph{holomorphic Euler-Poincar\'{e}
characteristic}, $e(S)$ the \emph{topological Euler number}, and
$K^2_S$ the \emph{self-intersection of the canonical divisor} (see e.g., \cite{Ba,BHPV,beau}).
Moreover, if $C$ is a smooth compact complex curve (Riemann surface),
then $g(C)$ will denote its \emph{genus}.

We shall also use a standard notation in group theory, hence we
denote by $\mathbb{Z}_n$ the \emph{cyclic group} of order $n$,
by $A_n$ the \emph{alternating group} on $n$ letters, by $S_n$ the
\emph{symmetric group} on $n$ letters, by $D_n$ the \emph{dihedral
group} of order $2n$, by $Q_8$ the \emph{group of quaternions}, by
$D_{p,q,r}$ a group with following presentation $\langle x,y \mid
x^p=y^q=1, \mbox{ }xyx^{-1}=y^r\rangle$ and $(r,q)=1$, by $\GL(2,q)$ the group
of invertible $2 \times 2$ matrices over the finite field with $q$
elements, which we denote by $\mathbb{F}_q$, and by $\SL(2,q)$ the subgroup
of $\GL(2,q)$ comprising the matrices with determinant $1$. With $Z(G)$ we shall
denote the center of a group $G$, moreover let $H \leq G$ be a
subgroup then the normalizer of $H$ in $G$ it will be denoted by
$N_G(H)$, while by $C_G(x)$ we denote the centralizer of $x \in
G$. In addition we shall write $x \sim_G y$ if $x,y \in G$ are
conjugate to each other.

%%%%%%%%%%%%%%%%%%%%%%%%%%%%%%%%%%%%%%%%%%%%%%%%%%%%%%%%%%%%%%%%%%%%%%%%%%%%%%%%%%%%%%%%%%%%%%%%%%%%%%%%%%%%%%%%%%%%%%%%%%%%%%
%--------------------------------------------------------------------------%
\section{Group theoretical preliminaries}
The study of surfaces isogenous to a product of curves is
strictly linked with the study of Galois coverings of Riemann surfaces
with fixed branching data. We collect in this section some
standard facts on coverings of Riemann surfaces from an algebraic
point of view following the notation of \cite{pol1} and \cite{catdd}.
%--------------------------------------------------------------------------%
\begin{defin}\label{genrvect} Let $G$ be a finite group and
let:
\[  0 \leq g', \ \ \ \ \ 2 \leq m_1 \leq \dots \leq m_r
\]
be integers. A \emph{generating vector} for $G$ of type
$(g'\mid m_1,\dots,m_r)$ is a $(2g'+r)-$tuple of elements of $G$:
\[ \mathcal{V}=\{a_1,b_1,\dots,a_{g'},b_{g'},c_1,\dots,c_r\}
\]
such that the following are satisfied:
\begin{enumerate}
    \item The set $\mathcal{V}$ generates $G$.
    \item Denoting by $| c |$ the order of $c$:
    \begin{description}
        \item[A] $| c_i | = m_i \ \ \forall 1\leq i \leq r$, or
        \item[B] there exist a permutation $\sigma \in S_r$ such that:
\[  |c_1| = m_{\sigma(1)},\dots | c_r | = m_{\sigma(r)},
\]
    \end{description}
    \item $\prod^{g'}_{i=1}[a_i,b_i]c_1 \dots  c_r=1$.
\end{enumerate}
If such a $\mathcal{V}$ exists we say that $G$ is $(g'\mid
m_1,\dots ,m_r)-$generated. We refer to $\mathbf{m}:=m_1,\dots,m_r$ as the \emph{branching
data}.
\\ Moreover if $g'=0$ a generating vector is said to be \emph{spherical}.
\end{defin}
%-----------------------------------------------------------------------------%
When we consider the definition with $(ii)$ \textbf{A} we shall
clearly speak of \emph{ordered vectors}. \emph{Unordered vectors} ($(ii)$
\textbf{B}) are needed only when we tackle the
problem of the moduli space, and so until the last section we
shall suppose that the generating vectors are all ordered.
%--------------------------------------------------------------------------%
We shall also use the notation, for example, $(g' \mid 2^4,3^2)$ to indicate
the tuple $(g' \mid 2,2,2,2,3,3)$.

We have the following reformulation of the Riemann Existence Theorem.
%--------------------------------------------------------------------------%
\begin{prop}\label{RiemEx} A finite group $G$ acts as a
group of automorphisms of some compact Riemann surface $C$ of
genus $g$ if and only if there exist integers $g' \geq 0$ and $m_r
\geq m_{r-1} \geq \dots \geq m_1\geq 2$ such that $G$ is $(g'\mid
m_1,\dots,m_r)-$generated for some generating vector
$\{a_1,b_1,\dots,a_{g'},b_{g'},c_1,\dots,c_r\}$, and the following
Riemann-Hurwitz relation holds:
\begin{equation}\label{RH} 2g-2=| G | (2g'-2 +
\sum^r_{i=1}(1-\frac{1}{m_i})).
\end{equation}
\end{prop}
If this is the case, then $g'$ is the genus of the quotient
Riemann surface $C':=C/G$ and the Galois covering $C \rightarrow C'$ is
branched in $r$ points $p_1,\dots,p_r$ with branching numbers
$m_1,\dots,m_r$ respectively. Moreover if $r=0$ the covering is said
to be \emph{unramified} or \emph{\'{e}tale}.

%--------------------------------------------------------------------------%
One introduces the following abstract group.
%--------------------------------------------------------------------------%
\begin{defin}\label{fuchs} Let us denote by $\Gamma = \Gamma(g' \mid
m_1,\dots,m_r)$ the abstract group with presentation:
\begin{multline*}  \Gamma(g' \mid m_1, \dots , m_r):=\langle \alpha_1,\beta_1, \ldots , \alpha_{g'},\beta_{g'},\gamma_{1},
\ldots , \gamma_{r} | \\
\gamma^{m_1}_{1}=\dots= \gamma^{m_r}_{r}=\prod_{k=1}^{g'} [\alpha_{k},\beta_{k}]\gamma_{1}
\cdot \ldots \cdot \gamma_{r} =1 \rangle.
\end{multline*}
We shall call this group an \emph{orbifold surface group}, following \cite{catdd}.
\end{defin}
%--------------------------------------------------------------------------%
Notice that other authors call this a \emph{Fuchsian type group} of signature $(g' \mid
m_1,\dots,m_r)$, see e.g., \cite{breuer}. Indeed, one can reinterpret the Riemann
Existence Theorem in terms of exact sequence of groups via
the Uniformization Theorem.
%--------------------------------------------------------------------------%
\begin{prop}\label{REII} A finite group $G$ acts as a group of
automorphisms on some compact Riemann surface $C$ of genus $g \geq
2$ if and only if there exists an exact sequence of groups:
\[ 1 \longrightarrow \Pi_g \longrightarrow \Gamma \stackrel{\theta}{\longrightarrow} G
\longrightarrow 1,
\]
where $\Gamma = \Gamma(g' \mid m_1,\dots,m_r)$ is an orbifold surface group, and $\Pi_g$ is the fundamental group
of $C$.
\end{prop}
Since $\Pi_g$ is torsion free, if follows that the order of each
generators  $\gamma_i \in \Gamma$ is the same as the one of $\theta(\gamma_i) \in G$ (see \cite[Lemma 3.6]{breuer}).
This is a reason to
give the following definition.
%--------------------------------------------------------------------------%
\begin{defin}\label{admepi} Let $G$ be a finite group. An
epimorphism $\theta: \Gamma = \Gamma(g' \mid m_1,\dots,m_r)
\rightarrow G$ is called \emph{admissible} if $\theta (\gamma_i)$ has order
$m_i$ for all $i$. If an admissible epimorphism exists, then clearly $G$
is $(g'\mid m_1,\dots,m_r)-$generated.
\end{defin}
%--------------------------------------------------------------------------%
One notices that the image of the generators of an orbifold surface group
under an admissible epimorphism are exactly the
generating vectors of Definition \ref{genrvect}.

Thanks to Propositions \ref{RiemEx} we are able to translate the geometrical
problem of finding Galois coverings $C \rightarrow C/G$ into the algebraic
problem of finding generating vectors for $G$ of type
$(g(C/G)\mid m_1,\dots ,m_r)$.
%--------------------------------------------------------------------------%
%--------------------------------------------------------------------------%
\section{Surfaces isogenous to a product with $p_g=q=2$ }
%--------------------------------------------------------------------------%
Recall the following definitions.
%------------------------------------------------------------------------------%
\begin{defin} A surface $S$ is said to be of Albanese general type, if $dim(\alpha (S))=2$, where $\alpha: S \rightarrow Alb(S)$
is the Albanese map.
\end{defin}
Notice that if $q(S)=2$ then $dim(Alb(S))=2$, thus the Albanese map is surjective if and only if $S$ is of Albanese of general type.  
%--------------------------------------------------------------------------%
\begin{defin}\label{def.isogenous} A surface $S$ is said to be \emph{isogenous to a (higher)
product of curves} if and only if, equivalently, either:
\begin{enumerate}
    \item $S$ admits a finite unramified covering, which is
    isomorphic to a product of curves of genera at least two, or
    \item $S$ is a quotient $S:=(C \times F)/G$ where $C$, $F$ are
    curves of genus at least two, and $G$ is a finite group acting
    freely on $C \times F$.
\end{enumerate}
\end{defin}
%--------------------------------------------------------------------------%
By Proposition 3.11 of \cite{cat00} the two properties
(\textit{i}) and (\textit{ii}) are equivalent. Using the same notation as in definition \ref{def.isogenous}, let
$S$ be a surface isogenous to a product, and $G^{\circ}:=G
\cap(Aut(C) \times Aut(F))$. Then $G^{\circ}$ acts on the two
factors $C$, $F$ and diagonally on the product $C \times
F$. If $G^{\circ}$ acts faithfully on both curves, we say that
$S= (C \times F)/G$ is a \emph{minimal
realization}. In \cite{cat00}
is also proven that any
surface isogenous to a product
admits a unique minimal realization. From now on we shall
work only with minimal realization.

By \cite[Lemma 3.8]{cat00} there are two cases: the
\emph{mixed}\index{Mixed Type} case where the action of $G$
exchanges the two factors, in this case $C$ and $F$ are
isomorphic and $G^{\circ} \neq G$; and the \emph{unmixed}\index{Unmixed Type} case, where
$G=G^{\circ}$ and therefore it acts diagonally.

We observe that a surface isogenous to a product of curves is of
general type, it is always minimal (see \cite{cat00} Remark
3.2), and its numerical invariants are explicitly given in terms of the genera of the curves and the order of the
group by the following proposition.
\begin{prop}[\cite{cat00}, Theorem 3.4] Let $S=(C \times F)/G$ be a surface isogenous to a product and denote by $d$ the order of $G$, then:
\begin{equation}\label{eq.euler.isot.fib}
e(S)=\frac{4(g(C)-1)(g(F)-1)}{d},
\end{equation}
\begin{equation}\label{eq.k2.isot.fib}
K^2_S=\frac{8(g(C)-1)(g(F)-1)}{d},
\end{equation}
\begin{equation}\label{eq.chi.isot.fib}
\chi(S)=\frac{(g(C)-1)(g(F)-1)}{d}.
\end{equation}
\end{prop}

From now on let $S$ be a surface of general type with
$p_g=q=2$ and not of Albanese general type.
According to \cite{cat00} we give the following definition.
%--------------------------------------------------------------------------%
\begin{defin} A surface isogenous to a product of unmixed type $S:=(C \times F)/G$
is said to be \emph{generalized hyperelliptic} if:
\begin{enumerate}
    \item the Galois covering $C \rightarrow C/G
    $ is unramified,
    \item the quotient curve $F/G$ is isomorphic to
    $\mathbb{P}^1$.
\end{enumerate}
\end{defin}
%--------------------------------------------------------------------------%
The following theorem gives a characterization of generalized hyperelliptic surfaces.
%--------------------------------------------------------------------------%
\begin{theo}[\cite{cat00}, Theorem 3.18]\label{Catatrediciotto} Let $S$ be a surface such
that:
\begin{enumerate}
    \item $K^2_S=8\chi (\mathcal{O}_S) > 0$
    \item $S$ has irregularity $q \geq 2$ and the Albanese map is
    a pencil.
\end{enumerate}
Then, letting $g$ be the genus of the Albanese fibres, we have:
$(g-1) \leq \frac{\chi (\mathcal{O}_S)}{q-1}$.
\\ A surface $S$ is generalized hyperelliptic if and only if
(i) and (ii) hold and $g=1+\frac{\chi (\mathcal{O}_S)}{q-1}$.
\\ In particular, every $S$ satisfying (i) and (ii) with $p_g=2q-2$
is generalized hyperelliptic where $G$ (from the
definition) is a group of automorphisms of the curve $F$ of genus
$2$ with $F/G \cong \mathbb{P}^1$.
\end{theo}
%--------------------------------------------------------------------------%
\begin{cor}[Zucconi, \cite{zucc} Prop. 4.2]\label{Zucco} If $S$ be a surface of general type
with $p_g=q=2$ and not of Albanese general type. Then $S$ is generalized hyperelliptic.
\end{cor}
%--------------------------------------------------------------------------%
\begin{rem}\label{main} We collect all the properties of
a surface $S$  of general type with $p_g=q=2$ not of Albanese
general type. Let $\alpha: S \rightarrow Alb(S)$ the Albanese map
and $B:=\alpha(S)$. Then:
\begin{enumerate}
    \item $S$ is isogenous to an unmixed product $(C \times F)/ G$.
    \item $K^2_S=8$.
    \item $g(F)=2$ and $F/G \cong \mathbb{P}^1$.
    \item $C \rightarrow C/G$ is unramified and $C/G \cong
    B$ has genus $2$.
    \item $| G | = (g(C)-1)(g(F)-1)= (g(C)-1)$.
\end{enumerate}
\end{rem}

%--------------------------------------------------------------------------%

To classify all the groups and the genera of smooth curves of
surfaces isogenous to a product with $p_g=q=2$ and not of Albanese
general type one can proceed as follows: first one classifies all
possible finite groups $G$ which induce a Galois covering $f\colon F
\rightarrow \mathbb{P}^1$ with $g(F)=2$, second one has to check
whether such groups $G$ induce an unramified Galois covering $g\colon
C \rightarrow B \cong C/G$, where the genus of $B$ is $2$ and the
genus of $C$ is determined by the Riemann-Hurwitz formula.

We notice that the action of $G$ on the
product $C \times F$ is always free, since the action on $C$ is free.
%--------------------------------------------------------------------------%
\begin{theo}\label{mainthone} Let $S$ be a complex surface of general type with
$p_g=q=2$  not of Albanese general type. Then $S=(C \times F)/G$ is
generalized hyperelliptic and assuming w.l.o.g. $g(F)=2$
the only possibilities for the genus of $C$, the group $G$ and the
branching data $\mathbf{m}$ for $F \rightarrow F/G \cong
\mathbb{P}^1$ are given by the entries in Table 1 of Theorem
\ref{zerozero} labelled with GH.
\end{theo}
%------------------------------------------------------------------------------%
%%%%%%%%%%%%%%%%%%%%%%%%%%%%%%%%%%%%%%%%%%%%%%%%%%%%%%%%%%%%%%%%%%%%%%%%%%%%%%%%%%%%%%%%%%%%%%%%
%%%%%%%%%%%%%%%%%%%%%%%%%%%%%%%%%%%%%%%%%%%%%%%%%%%%%%%%%%%%%%%%%%%%%%%%%%%%%%%%%%%%%%%%%%%%%%%%
%
%%%%%%%%%%%%%%%%%%%%%%%%%%%%%%%%%%%%%%%%%%%%%%%%%%%%%%%%%%%%%%%%%%%%%%%%%%%%%%%%%%%%%%%%%%%%%%%%
%%%%%%%%%%%%%%%%%%%%%%%%%%%%%%%%%%%%%%%%%%%%%%%%%%%%%%%%%%%%%%%%%%%%%%%%%%%%%%%%%%%%%%%%%%%%%%%%
\begin{proof}
We are in the hypothesis of Corollary \ref{Zucco}, thus $S=(C \times F) /G$, where $F$, $C$ and $G$ have the
property indicated in Remark \ref{main}.

The classification of the automorphisms groups of a Riemann surface
$F$ of genus $2$ was given by Bolza in \cite{bol}, moreover the
classification of all the groups $G$ acting effectively as a group
of automorphisms of $F$ such that the quotient $F/G$ is isomorphic
to $\mathbb{P}^1$ is given in \cite{zucc} or \cite{Broug}. We give
a full proof of the classification of the latter groups, since we
are interested in obtaining a complete information including also the
branching data.

By the Riemann-Hurwitz formula (\ref{RH}) applied to $F$ we obtain:
\begin{equation}\label{HR1} 2g(F)-2=|G| \big(-2 + \sum^r_{i=1}(1-\frac{1}{m_i})\big),
\end{equation}
remembering that $g(F)=2$ we get:
\begin{equation}\label{HR4}
2= |G |\big(-2 + \sum^r_{i=1}(1-\frac{1}{m_i})\big)
\end{equation}
which yields:
\[ |G| (\frac{r}{2}-2) \leq 2 \leq |G|(r-2),
\]
and since $| G | \geq 2$ we have $3 \leq r \leq 6 $.
\\ We examine all the cases proceeding as follows: for each $r$, using the fact that $2 \leq m_1 \leq \dots \leq m_r$, and by (\ref{HR4}), we can bound the order of $G$ from above by a rational function of $m_1$  and from below by $m_1$, since $m_1$ divides $|G|$, and we analyze
case by case. As soon as $m_1$ gives a void condition, we repeat the same analysis using $m_2$, and so on for all $m_i$'s. In the first case we shall perform a full calculation as an example.
%%%%%%%%%%%%%%%%%%%%%%%%%%%%%%%%%%%%%%%%%%%%%%%%%%%%%%%%%%%%%

 \textbf{Case} $\mathbf{r=6}$.
%%%%%%%%%%%%%%%%%%%%%%%%%%%%%%%%%%%%%%%%%%%%%%%%%%%%%%%%%%%%
\\ In this case by (\ref{HR4}) we have :
\[ 2= | G | \big(-2 + \sum^6_1(1-\frac{1}{m_i})\big) \geq | G | \big(-2 + 6 (1-\frac{1}{m_1})\big),
\]
which yields
\[ m_1 \leq | G | \leq \frac{m_1}{2m_1-3},
\]
then $m_1=2$ and $|G| = 2$. Therefore by equation
(\ref{eq.chi.isot.fib}) $g(C)=3$, and since $m_i$ divides $| G
|$ for all $i=1,\dots,6$, we have $m_i=2$ for all $i=1,\dots,6$.
Then $G = \mathbb{Z}_2$, since it is $(0 \mid 2^6)-$generated.
We recover the first case in Table 1, i.e., $g(F)=2$ $g(C)=3$ and
$\mathbf{m}=(2^6)$.
\\ Notice that to fully recover this first case we still have to prove that $\mathbb{Z}_2$ induces an unramified covering $g\colon C \rightarrow B \cong C/\mathbb{Z}_2$, where the genus of $B$
is $2$, but this is obvious. From now on in order to avoid many repetitions we investigate the branching data and the order of the groups and it will be clear which case in Table 1 is recovered. Moreover we shall prove at the end that all the groups, we have found, induce an unramified covering  $C\rightarrow C/G $ with quotient a curve of genus $2$.
%%%%%%%%%%%%%%%%%%%%%%%%%%%%%%%%%%%%%%%%%%%%%%%%%%%%%%%%%%%%%%

 \textbf{Case} $\mathbf{r=5}$.
%%%%%%%%%%%%%%%%%%%%%%%%%%%%%%%%%%%%%%%%%%%%%%%%%%%%%%%%%%%%%
\\Proceeding as in the previous case, we have:
\[ 2 \leq | G | \leq 4.
\]
If $| G | =2$ then $m_i=2$ for all $i=1,\dots,5$ which yields
a contradiction to (\ref{HR4}).
\\ If $| G | = 3$ then
$m_i=3$ for all $i=1,\dots,5$ and again we have a contradiction.
\\ If $|G| =4$ then $m_i=2$ for all $i=1,\dots,5$%, because of (\ref{HR4})
. Since the elements of order $2$ generate the group we have $G =
(\mathbb{Z}_2)^2$. Indeed we have $c_1, \dots c_5, \in
(\mathbb{Z}_2)^2\setminus \{(0,0)\}$ such that
$\sum^5_{i=1}c_i=0$, for example take $c_1$, $c_2$ and $c_3$ all
different from each other, and $c_1=c_4=c_5$.
%%%%%%%%%%%%%%%%%%%%%%%%%%%%%%%%%%%%%%%%%%%%%%%%%%%%%%%%%%%%%

 \textbf{Case} $\mathbf{r=4}$.
%%%%%%%%%%%%%%%%%%%%%%%%%%%%%%%%%%%%%%%%%%%%%%%%%%%%%%%%%%%%%%
\\ In this case we have:
\[ m_1 \leq | G | \leq \frac{2m_1}{2m_1-4},
\]
then $m_1 \leq 3$.
\\ If $m_1=3$, then $| G |=3$ and $m_i=3$ for all $i=1,\dots 4$.
Clearly $G=\mathbb{Z}_3$ which is $(0\mid 3^4)-$generated,
consider for example $c_1=1$, $c_2=1$, $c_3=2$, $c_4=2$.
\\ Now suppose $m_1=2$, this gives no upper bound for the order of $G$, therefore looking at the possible values of $m_2$, we have:
\[ \mbox{l.c.m.}(2, m_2) \leq | G | \leq \frac{4m_2}{3m_2-6}.
\]
We can exclude the cases with $m_2 \geq 3$, so $m_2=2$.
\\ If we proceed further and look at the values of $m_3$ once $m_1=m_2=2$, we have:
\[ \mbox{l.c.m.}(2, m_3) \leq | G | \leq \frac{2m_3}{m_3-2},
\]
so that $m_3 \leq 4$.
\\ If $m_3=4$, then $| G | =4$, $m_4=4$ and $G=\mathbb{Z}_4$, which
is $(0 \mid 2^2,4^2)-$generated, take for example $c_1=2$,
$c_2=2$, $c_3=1$, $c_4=3$.
\\ If $m_3=3$, then $| G |
=6$ and we have $m_4=3$. We have two possibilities either
$G=\mathbb{Z}_6$ or $G=S_3$, and both cases occur, since both
groups are $(0 \mid 2,2,3,3)-$generated. For the first case
consider for example $c_1=3$, $c_2=3$, $c_3=2$, and $c_4=4$, while
for the latter one $c_1=(1,2)$, $c_2=(2,3)$, $c_3=(1,3,2)$,
$c_4=(1,3,2)$.
\\ Let us consider the case $m_3=2$, then we have to look at the possible values of $m_4$, since:
\[ \mbox{l.c.m.}(2,m_4) \leq | G |=\frac{4m_4}{m_4-2}
\]
then only possibilities are the following:
\\ If $m_4=6$, then $| G | = 6$ and this case is impossible. Indeed $G$ cannot be $S_3$, since $S_3$ has no
element of order $6$. In addition let $c_1,\dots,c_4$ be the
generators of order $m_1,\dots,m_4$, then we must have $c_1 + c_2
+ c_3 + c_4=0$, then it cannot be $\mathbb{Z}_6$ since the only
element of order two is $3$ and $3$ plus an element of order six
is never $0$.
\\ If $m_4=4$ then $| G | = 8$ then $G= D_4$, for example with the following generators:
$c_1=y$, $c_2=yx$, $c_3=x^2$, $c_4=x$, where $y$ is a reflection
and $x$ a rotation. The group $G$ cannot be $\mathbb{Z}_2 \times
\mathbb{Z}_4$ or $Q_8$, since the conditions $c_1 + c_2 + c_3 +
c_4=0$, respectively $c_1\cdot c_2 \cdot c_3 \cdot c_4=1$ are not
satisfied. $G$ cannot be $\mathbb{Z}_8$, because it is not $(2,4)$
generated. $G$ cannot be $(\mathbb{Z}_2)^3$ since it does not have
any element of order $4$.
\\ If $m_4=3$ then $| G | = 12$, $G$ cannot be
$\mathbb{Z}_{12}$ or $\mathbb{Z}_2 \times \mathbb{Z}_6$ because of
the condition $c_1+\dots+c_4=0$, $G$ cannot be
$D_{3,4,-1}:=\langle x,y \mid x^4=y^3=1, \mbox{ } xyx^{-1}=y^{-1}
\rangle$ because $D_{3,4,-1}$ has only one element of order $2$.
The two remaining cases are $D_6$ or $A_4$. If $G$ has four
$3-$Sylow subgroups then $G=A_4$, impossible since the elements of
order $2$ are in the Klein subgroup while $c_4$ is not. In the
other case $\mathbb{Z}_3 \subset G$ is normal, and there is an
element of order $6$. We recover the case $G= D_6$ with generating
vector: $c_1=y$, $c_2=yx$, $c_3=x^3$, $c_4=x^2$.
%%%%%%%%%%%%%%%%%%%%%%%%%%%%%%%%%%%%%%%%%%%%%%%%%%%%%%%%%%%%%%%%%%%%%%%%%%%%%%%%%%%%%%%

\textbf{Case} $\mathbf{r=3}$.
%%%%%%%%%%%%%%%%%%%%%%%%%%%%%%%%%%%%%%%%%%%%%%%%%%%%%%%%%%%%%%%%%%%%%%%%%%%%%%%%%%%%%%%
\\ This case is much more involved
than the previous ones. We have
\[ m_1 \leq | G | \leq \frac{2m_1}{m_1-3},
\]
then after a short calculation one sees that $m_1 \leq 5$.
%%%%%%%%%%%%%%%%%%%%%%%%%%%%%%%%%%%%%%%%%%%%%%%%%%%%%%%%%%%%%%%%%%%%%%%%%%%%%%%%%%%%%%%%

If $\mathbf{m_1=5}$ then $| G |=5$ and $m_2=m_3=5$. The only
possibility is $G= \mathbb{Z}_5$ which is $(0\mid 5^3)-$generated,
for example consider $c_1=1$, $c_2=2$ and $c_2=2$.
%%%%%%%%%%%%%%%%%%%%%%%%%%%%%%%%%%%%%%%%%%%%%%%%%%%%%%%%%%%%%%%%%%%%%%%%%%%%%%%%%%%%%%%%

If $\mathbf{m_1=4}$ one has that $| G |=8$ and $m_2=m_3=4$, and
the only group of order $8$ which is $(0 \mid 4,4,4)-$generated is
$G=Q_8$, consider for example $c_1=i$, $c_2=j$ and $c_3=-k$.
Notice that the other groups of order $8$ containing an element of
order $4$ are $\mathbb{Z}_8$, $\mathbb{Z}_2 \times \mathbb{Z}_4$, and $D_4$.
The elements of order $4$ in $\mathbb{Z}_8$ and in $D_4$ form  proper
subgroups, so these cases are excluded. In $\mathbb{Z}_2 \times
\mathbb{Z}_4$ there are four elements of order $4$, the sum of any two
of them is an element of order at most $2$, so condition $c_1 +
c_2 + c_3=0$ cannot be satisfied.
%%%%%%%%%%%%%%%%%%%%%%%%%%%%%%%%%%%%%%%%%%%%%%%%%%%%%%%%%%%%%%%%%%%%%%%%%%%%%%%%%%%%%%%%%%%%%%

If $\mathbf{m_1=3}$ we have to look at all possible values of
$m_2$, since:
\[ \mbox{l.c.m.}(3, m_2) \leq | G | \leq  \frac{3m_2}{m_2-3}
\]
then only possibilities are the following:
\\ If $m_2=6$, then $| G |
=6$ and $m_3=6$. The only possibility is $G=\mathbb{Z}_6$,
consider for example as generators $c_1=4$, $c_2=1$ and $c_3=1$.
Notice that $S_3$ has no elements of order $6$.
\\ If $m_2=4$, then $| G |=12$ and $m_3=4$. In
this case the only group of order $12$ which can be generated by
elements $c_1,c_2,c_3$ of order $3,4,4$ and such that these
elements satisfy
\begin{equation}\label{somazero}
c_1 \cdot c_2 \cdot c_3=1 \textrm{ (or additively } c_1 + c_2 +
c_3=0),
\end{equation}
is $D_{4,3,-1}$, choose for example $c_1=y$, $c_2=xy$ and
$c_3=x^3$, where the notation is the one given above. Notice that
all the other groups either do not have an element of order $4$ or
it is $\mathbb{Z}_{12}$, which fails condition (\ref{somazero}).
\\ In the case $m_2=3$ we have to look at the possible values of $m_3$, since:
\[ \mbox{l.c.m.}(3,m_3) \leq | G | = \frac{6m_3}{m_3-3},
\]
then the only possibilities are the following:
\\ If $m_3=9$, then $| G | =9$ and $G$ could be only $\mathbb{Z}_{9}$, but condition (\ref{somazero}) is not satisfied by elements of order $3,3,9$, which excludes this case.
\\ If $m_3=6$, then $| G | =12$. Also this case has to be
excluded because: $A_4$ does not have an element of order $6$,
while for $D_6$ and $(\mathbb{Z}_2)^2 \times \mathbb{Z}_3$ the
elements of order $3$ and $6$ cannot generate, in the end for the
last two groups (\ref{somazero}) fails.
\\ If $m_3=5$ then $| G | =15$ and we have only
$\mathbb{Z}_{15}$, but its generating elements of order $3$ and
$5$ do not satisfy (\ref{somazero}).
\\ If  $m_3=4$ then $| G | =24$. Here the
number of groups involved or their order can be considerably large.
In order to avoid many repetitions if these numbers are
excessively large, where indicated, we use a computer program in
GAP4 (see Appendix for an explanation of the program),
to check the corresponding cases. Indeed the computer
shows that among the $15$ groups of order $24$ the only one which
can be $(0 \mid 3,3,4)-$generated is \texttt{SmallGroup(24,3)},
which corresponds to $G=\SL(2,\mathbb{F}_3)$. This exhausts all the cases with
$m_1=3$.
%%%%%%%%%%%%%%%%%%%%%%%%%%%%%%%%%%%%%%%%%%%%%%%%%%%%%%%%%%%%%%%%%%%%%%%%%%%%%%%%%%%%%%%%%%

If $\mathbf{m_1=2}$ we look at all possible values of $m_2$ since:
\[ \mbox{l.c.m.}(2,m_2) \leq  | G | \leq \frac{4m_2}{m_2-4},
\]
then only possibilities are the following:
\\ If $m_2=8$, then $| G | = 8$
and $m_3=8$ which yields $G=\mathbb{Z}_8$, choose for example as
generating vector $c_1=4$, $c_2=5$, $c_3=7$.
\\ If $m_2=6$, then $| G | =12$ and $m_3=6$ which yields the case:
$G= \mathbb{Z}_2 \times \mathbb{Z}_6$, choose for example as
generating vector $c_1=(1,3)$, $c_2=(1,2)$, $c_3=(0,1)$. Notice
that it cannot be $\mathbb{Z}_{12}$ because of (\ref{somazero}),
$G$ cannot be $A_4$ because it does not have any element of order
$6$. Moreover $G$ cannot be $D_6$, since to generate it one needs
a reflection $y$ but the condition $c_2 c_3=y$ can never hold
since the only elements with order $6$ are $x$ and $x^5$, with $x$
rotation. Finally $D_{4,3,-1}$ is impossible because the two
elements of order $6$ and the element of order $2$ do not satisfy
(\ref{somazero}).
\\ If $m_2=5$ then $| G |  \leq 20$, looking at the branching data the only two possible cases
are $| G | =20, 10$.
\\ If $| G | =20$, then $m_3=5$.
Among the $5$ groups of order $20$ a computer computation shows
that none of them are $( 0 \mid 2,5,5)-$ generated.
\\ If $| G | =10$ then $m_3=10$, which gives $G=\mathbb{Z}_{10}$ for example if
$c_1=x5$, $c_2=4$, $c_3=1$. $G$ cannot be $D_5$
since it has no element of order $10$.
\\ If $m_2 \leq 4$ one has to look at all possible values of $m_3$.
\\ Let $m_2=4$, since:
\[ \mbox{l.c.m.}(4,m_3) \leq | G | = \frac{8m_3}{m_3-4}
\]
then the only possibilities are the following:
\\ If $m_3=12$, then $| G | =12$  and $G$ could be only  $\mathbb{Z}_{12}$, but condition (\ref{somazero}) cannot be satisfied by elements of order $2$, $4$, $12$, therefore this case is excluded.
\\ If $m_3=8$, then $| G | =16$,
among the $14$ groups of order $16$ a GAP4 computation shows that
only \texttt{SmallGroup(16,8)} (i.e., $G =D_{2,8,3}$) can be $(0 \mid 2,4,8)-$generated.
\\ If $m_3=6$, then
$| G | = 24$ and the only possibility for $G$ is \texttt{SmallGroup(24,8)}, which is
$G= \mathbb{Z}_2 \ltimes
((\mathbb{Z}_2)^2 \times \mathbb{Z}_3)$. This case was also accomplished using GAP4.
\\ One sees that case $m_3=5$ is impossible, here again it is needed a computational fact: none of the $14$ groups of
order $40$ can be $(0 \mid 2,4,5)$ generated.
\\ We now consider the case $m_2=3$. We look at all possible values of $m_3$, since:
\[ \mbox{l.c.m.}(6,m_3) \leq | G | = \frac{12m_3}{m_3-6}
\]
then only possibilities are the following.
\\ If $m_3=18$, then $| G | =18$  and $G = \mathbb{Z}_{18}$, but condition (\ref{somazero})
cannot be satisfied by elements of order $2,3,18$.
\\ If $m_3=12$, then $| G | =24$, a computer calculation shows that among the $15$ groups of order $24$
none of them are $(0 \mid 2,3,12)-$generated.
\\ We can also exclude the case $m_3=10$ because none of the groups
of order $30$ is $(0 \mid 2,3,10)-$generated.
\\ One can exclude the case $m_3=9$, because among the $14$ groups of order $36$ a
computer computation shows that none of them are $(0 \mid
2,3,9)-$generated.
\\ Case $m_3=7$ is
also excluded, though there are $15$ groups of order $84$,
a computer computation shows that none of them can be $(0 \mid
2,3,7)-$generated.
\\ The remaining case is  $m_3=8$ which gives $| G | = 48$.
A computer computation shows that among the $52$
groups of order $48$ only  \texttt{SmallGroup(48,29)} (i.e., $G= \GL(2, \mathbb{F}_3)$) satisfies all the necessary conditions. This
recovers the last case of Table 1.
%%%%%%%%%%%%%%%%%%%%%%%%%%%%%%%%%%%%%%%%%%%%%%%%%%%%%%%%%%%%%%%%%%%%%%%%%%%%%%%%%%%%%%%%%

Now we have to see whether for each possible group $G$ there is a
surjective homomorphism:
\[ \Gamma(2 \mid -) \twoheadrightarrow G.
\]
Indeed this is true in all the cases, more precisely one notices
that all the possible groups $G$ can be two generated, call the generators $x$ and
$y$. Then we have the following epimorphism from
$\Gamma(2 \mid -)$ to $G$:
\[ \alpha_1 \mapsto x, \ \ \ \beta_1 \mapsto 1, \ \ \ \alpha_2 \mapsto y,  \ \ \ \beta_2 \mapsto 1.
\]
\end{proof}
%%%%%%%%%%%%%%%%%%%%%%%%%%%%%%%%%%%%%%%%%%%%%%%%%%%%%%%%%%%%%%%%%%%%%%%%%%%%%%%%%%%%%%%%%%%%%%%%
%%%%%%%%%%%%%%%%%%%%%%%%%%%%%%%%%%%%%%%%%%%%%%%%%%%%%%%%%%%%%%%%%%%%%%%%%%%%%%%%%%%%%%%%%%%%%%%%
%
%%%%%%%%%%%%%%%%%%%%%%%%%%%%%%%%%%%%%%%%%%%%%%%%%%%%%%%%%%%%%%%%%%%%%%%%%%%%%%%%%%%%%%%%%%%%%%%%
%%%%%%%%%%%%%%%%%%%%%%%%%%%%%%%%%%%%%%%%%%%%%%%%%%%%%%%%%%%%%%%%%%%%%%%%%%%%%%%%%%%%%%%%%%%%%%%%

We have analyzed the case when the image of the Albanese map is a
curve. Now we want to see whether there are surfaces isogenous to a product of unmixed type with $p_g=q=2$ and of Albanese general type.

Let $S \rightarrow B$ an isotrivial fibration with general fibre $F$ and $g(F)\geq 2$, then a monodromy argument shows that $S$ is birational to the quotient of a product of two curves $(C \times F)/G$ modulo the diagonal action of a finite group $G$ and $C/G \cong B$ (see e.g., \cite[\S 1.1]{serr}). Moreover by \cite[Proposition 2.2]{serr} we have:
\begin{equation}\label{eq.irregularity} q(S)=g(C/G)+g(F/G).
\end{equation}
This applies in particular if $S$ is a surface isogenous to a product of unmixed type. If $q(S)=2$,  there are two cases: either (w.l.o.g) $g(C/G)=2$ and $g(F/G)=0$, or $g(C/G)=1$ and $g(F/G)=1$. 
By \cite[Proposition 4.3]{zucc} the first case is completely solved.
%\begin{prop}
%Let $C$ and $F$ be two smooth curves and let $G$ be a non-trivial
%finite group with two injections: $G \hookrightarrow Aut(C)$, $G
%\hookrightarrow Aut(F)$. Suppose $g(F)=2$, $F/G =\mathbb{P}^1$ and
%$\pi\colon C \rightarrow C/G$ is an \'{e}tale morphism where
%$g(C/G)=2$. Then the quotient $S=(C \times F)/G$ by the diagonal
%action is a minimal surface of general type with $p_g(S)=q(S)=2$
%and non-surjective Albanese morphism.
%\end{prop}
%----------------------------------------------------------------------------------------%

For the second case we search for surfaces
isogenous to an unmixed product $S=(C \times F)/G$ such that $C/G$
and $F/G$ are both elliptic curves and $\chi(S)=1$. We need the following two
results to simplify our search.
\begin{lem}\cite[Lemma 2.3, Corollary 2.4]{zucc}\label{duecinque} Let $S$ be surface of Albanese general type
with $p_g=q=2$. Let $\phi\colon S \rightarrow B$ be a fibration of
curves of genus $g$. If the genus of $B$ is $b > 0$, then $b=1$
and $ 2 \leq g \leq 5$.
\end{lem}
%And this fact about generating vectors.
\begin{lem}\label{noabel}
If $G$ is an abelian group and $G$ is $(g' \mid m_1,\dots,m_r)$-generated, then $r \neq 1$.
\end{lem}
\begin{proof}
Suppose $G$ abelian and $r=1$. Then the relation
$\Pi^{g'}_{i=1}[a_i,b_i]c_1=1$ yields $c_1=1$.
\end{proof}

In case the Albanese map is not surjective we have that $S$ is of
generalized hyperelliptic type, hence one of the two coverings is
\'{e}tale, and we have always a free action of $G$ on the product
$C \times F$. In case the Albanese map is surjective we do not
have an \'{e}tale covering, so we also have to check whether the
action of $G$ on the product of the two curves is free or not.
%----------------------------------------------------------------------%
\begin{rem}\label{singcong} Let  $\mathcal{V}_1:=(a_{1,1},b_{1,1},c_{1,1},
\dots ,c_{1,r_1})$ and $\mathcal{V}_2:=(a_{2,1},b_{2,1},c_{2,1},
\dots ,c_{2,r_2})$ be generating vectors of type $(1 \mid m_{i,1},
\dots ,m_{i,r_i})$ for $i=1,2$ respectively for a finite group $G$. Then the cyclic
groups $\left\langle c_{1,1}\right\rangle,\dots,\left\langle
c_{1,r_1}\right\rangle$ and their conjugates provide the non-trivial
stabilizers for the action of $G$ on $C$, whereas $\left\langle
c_{2,1}\right\rangle,\dots,\left\langle c_{2,r_2}\right\rangle$ and
their conjugates provide the non-trivial stabilizers for the
action of $G$ on $F$. The singularities of $(C \times F)/G$ arise
from the points of $C \times F$ with nontrivial stabilizer, since
the action of $G$ on $C \times F$ is diagonal, it follows that the
set $\mathcal{S}$ of all nontrivial stabilizer for the action of
$G$ on $C \times F$ is given by $\Sigma(\mathcal{V}_1) \cap
\Sigma(\mathcal{V}_2)$ where
\[ \Sigma(\mathcal{V}_i):= \bigcup_{h \in G} \bigcup^{\infty}_{j=0} \bigcup^{r_i}_{k=1} h \cdot c^j_{i,k} \cdot
h^{-1}.
\]
It is clear that if we want $(C \times F)/G$ to be smooth we have to require
$\mathcal{S}=\{1_G\}$. If this is the case we shall say that $\mathcal{V}_1$ and $\mathcal{V}_2$ are \emph{disjoint}.
\end{rem}
%--------------------------------------------------------------------%
\begin{theo} \label{unmixcase}
Let $S$ be a complex surface with $p_g=q=2$ of Albanese general
type and isogenous to a product of unmixed type. Then $S$ is
minimal of general type and the only possibilities for the genera of
the two curves $C$, $F$, the group $G$ and the branching data
$\mathbf{m}$ respectively for $F \rightarrow F/G$ and $C
\rightarrow C/G$ are given by the entries in Table 1 labelled with
UnMix.
\end{theo}
%---------------------------------------------------------------------%
%%%%%%%%%%%%%%%%%%%%%%%%%%%%%%%%%%%%%%%%%%%%%%%%%%%%%%%%%%%%%%%%%%%%%%%%%%%%%%%%%%%%%%%%%%%%%%%%
%%%%%%%%%%%%%%%%%%%%%%%%%%%%%%%%%%%%%%%%%%%%%%%%%%%%%%%%%%%%%%%%%%%%%%%%%%%%%%%%%%%%%%%%%%%%%%%%
%
%%%%%%%%%%%%%%%%%%%%%%%%%%%%%%%%%%%%%%%%%%%%%%%%%%%%%%%%%%%%%%%%%%%%%%%%%%%%%%%%%%%%%%%%%%%%%%%%
%%%%%%%%%%%%%%%%%%%%%%%%%%%%%%%%%%%%%%%%%%%%%%%%%%%%%%%%%%%%%%%%%%%%%%%%%%%%%%%%%%%%%%%%%%%%%%%%
\begin{proof}
By Lemma \ref{duecinque} we have, that $2 \leq g(F) \leq 5$ and, w.l.o.g., $g(F) \leq g(C)$. Moreover we have $g(C/G)=g(F/G)=1$. We analyze case by case:

$\mathbf{g(F)=2.}$
%%%%%%%%%%%%%%%%%%%%%%%%%%%%%%%%%%%%%%%%%%%%%%%%%%%%%%%%%%%%%%%%%%%%%%%%%%%%%%%%%%%%%%
\\ From
\begin{equation}\label{HR5}
2g(F)-2 = | G | \sum^r_{i=1}(1-\frac{1}{m_i})
\end{equation}
and $\sum^r_{i=1}(1-\frac{1}{m_i}) \geq \frac{1}{2}$  we have:
\[ 2 \leq | G | \leq 4,
\]
which yields:
\\ $|G|=4$ if and only if  $r=1$ and $m_1=2$;
\\ $|G|=3$ if and only if  $r=1$ and $m_1=3$;
\\ $|G|=2$ if and only if $m_1=m_2=2$.
\\ The first two cases contradict Lemma \ref{noabel}, the third one is
also impossible. First notice that, from equation
(\ref{eq.chi.isot.fib}), $g(C)=3$, and for $F$ and $C$ we have
respectively the following branching data: $(2,2)$ and
$(2,2,2,2)$. It follows that we do not have any free action of
$\mathbb{Z}_2$ on $C \times F$.

$\mathbf{g(F)=3.}$
%%%%%%%%%%%%%%%%%%%%%%%%%%%%%%%%%%%%%%%%%%%%%%%%%%%%%%%%%%%%%%%%%%%%%%%%%%%%%%%%%%%%%
\\ From equation $(\ref{HR5})$ we have:
\[ 2 \leq | G | \leq 8,
\]
moreover $2$ divides $| G |$ by equation (\ref{eq.chi.isot.fib}). Then we have to analyze the cases: $| G
| = 8,6,4,2$.
\\ If $| G | =8$ then $g(C)=5$, and by Riemann-Hurwitz the branching data for $F$ and $C$ are respectively $(2)$ and $(2,2)$.
By Lemma \ref{noabel} $G$ is not abelian and since it is $(1 \mid 2^2)-$generated it must be $D_4$. Indeed $G$ cannot be $Q_8$, because $Q_8$ is not $(1 \mid 2^2)-$generated, since the only element of order $2$ is $-1$.
One sees that $D_4$ acts on $C \times F$ freely, hence this case occurs. Indeed one can choose the following generating vectors:
\[ a_{1,1}=x \ \ b_{1,1}=y \ \ c_{1,1}=x^2;
\]
\[ a_{2,1}=x \ \ b_{2,1}=y \ \ c_{2,1}=x^2y \ \ c_{2,2}=y,
\]
where $x$ is a rotation and $y$ is a reflection. Then $\{x^2\}$ is one conjugacy class and $\{y, x^2y\}$ is another, hence the two vectors are disjoint.
\\ If $| G | =6$ then $g(C)=4$. The Riemann-Hurwitz formula yields $(3)$ as branching data for $F$ and
$(2,2)$ for $C$, which yield $G = S_3$. One sees that any pair of generating vectors is disjoint since the conjugacy classes of elements of order $2$ and $3$ are disjoint, thus $S_3$ acts on $C \times F$ freely, and this case occurs.
\\ If $| G | =4$ then $g(C)=3$. In this case the branching data of $F$ and $C$ are respectively $(2,2)$ and $(2,2)$. If $G=\mathbb{Z}_4$ the
action cannot be free, since there is only one element of order
two. For $\mathbb{Z}_2 \times \mathbb{Z}_2$ there is a free
action. We first notice, since $G$ is abelian,
$2c_1=2c_2=c_1+c_2=0$ and $2c'_1=2c'_2=c'_1+c'_2=0$, then we can
choose $c=c_1=c_2$ and $c'=c'_1=c'_2$. If we choose, for example,
$c=(1,1)$ and $c'=(1,0)$ we see that this case occurs.
\\ We observe that the case $| G | =2$ leads to a contradiction to $g(F) \leq g(C)$.

$\mathbf{g(F)=4.}$
%%%%%%%%%%%%%%%%%%%%%%%%%%%%%%%%%%%%%%%%%%%%%%%%%%%%%%%%%%%%%%%%%%%%%%%%%%%%%%%%%%5
\\From equation $(\ref{HR5})$ we have:
\[ 3 \leq | G | \leq 12,
\]
moreover $3$ must divide $| G |$ by equation (\ref{eq.chi.isot.fib}). Furthermore since we assumed $g(F) \leq g(C)$ the only remaining cases
are $| G | =12,9$.
\\ If $| G | =12$ then $g(C)=5$. The branching data of $F$ and $C$ are respectively $(2)$
and $(3)$. There is no non-abelian group of order  $12$ which is
simultaneously  $(1 \mid 2)$ and $(1 \mid 3)-$generated. To see
this one notices that the derived subgroups of $D_6$ and
$D_{3,4,-1}$ are both isomorphic to $\mathbb{Z}_3$, hence in both
cases there are no commutators of order $2$, therefore the two
groups cannot be $(1 \mid 2)-$generated. Moreover $A_4$ is not $(1
\mid 3)-$generated because its derived subgroup is $\mathbb{Z}_2
\times \mathbb{Z}_2$, hence there are no commutators of order $3$.
\\ If $| G | =9$ then $g(C)=4$. We see that the branching data for $F$ is $(3)$, since all
the groups of order $9$ are abelian, this case does not occur.

%%%%%%%%%%%%%%%%%%%%%%%%%%%%%%%%%%%%%%%%%%%%%%%%%%%%%%%%%%%%%%%%%%%%%%%%%%%%%%%%%%%
$\mathbf{g(F)=5.}$
\\ From equation $(\ref{HR5})$ we have:
\[ 4 \leq | G | \leq 16,
\]
moreover $4$ must divide $| G |$, and since $g(F) \leq g(C)$ the only case remaining is
$| G | =16$.
\\ If $ | G |= 16$ then $g(C)=5$, and the branching data for $F$ and $C$ are $(2)$
and $(2)$. Looking at the table in \cite{canpol} one sees that
among the $14$ groups of order $16$ only $\mathbb{Z}_4 \ltimes
(\mathbb{Z}_2)^2$,  $D_{4,4,-1}$ and $D_{2,8,5}$ are $(1 \mid
2)-$generated. Moreover with a computer computation using the
program of the Appendix, one sees that in these cases
the action of the groups on the product $C \times F$ cannot be
free. Therefore this case does not occur.
\end{proof}
%%%%%%%%%%%%%%%%%%%%%%%%%%%%%%%%%%%%%%%%%%%%%%%%%%%%%%%%%%%%%%%%%%%%%%%%%%%%%%%%%%%%%%%%%%%%%%%%
%%%%%%%%%%%%%%%%%%%%%%%%%%%%%%%%%%%%%%%%%%%%%%%%%%%%%%%%%%%%%%%%%%%%%%%%%%%%%%%%%%%%%%%%%%%%%%%%
%
%%%%%%%%%%%%%%%%%%%%%%%%%%%%%%%%%%%%%%%%%%%%%%%%%%%%%%%%%%%%%%%%%%%%%%%%%%%%%%%%%%%%%%%%%%%%%%%%
%%%%%%%%%%%%%%%%%%%%%%%%%%%%%%%%%%%%%%%%%%%%%%%%%%%%%%%%%%%%%%%%%%%%%%%%%%%%%%%%%%%%%%%%%%%%%%%%

Now we study the mixed case.
%--------------------------------------------------------------------------------------------%
\begin{theo}\cite[Proposition 3.16]{cat00} Assume that $G^{\circ}$
is a finite group satisfying the following properties:
\begin{enumerate}\renewcommand{\theenumi}{\it \roman{enumi}}
    \item $G^{\circ}$ acts faithfully on a smooth curve $C$ of genus
    $g(C) \geq 2$,
    \item There is a non split extension:
    \begin{equation}\label{split} 1 \longrightarrow G^{\circ} \longrightarrow G \longrightarrow
    \mathbb{Z}_2 \longrightarrow 1.
    \end{equation}
Let us fix a lift $\tau ' \in G$ of the generator of
$\mathbb{Z}_2$. Conjugation by $\tau '$ defines an element
$[\varphi]$ of order $\leq 2$ in $Out(G^{\circ})$.
\end{enumerate}
Let us choose a representative $\varphi \in Aut(G^{\circ})$ and
let $\tau \in G^{\circ}$ be such that $\varphi^2$ is equal to
conjugation by $\tau$. Denote by $\Sigma_C$ the set of elements in
$G^{\circ}$ fixing some point on $C$ and assume that both the
following conditions are satisfied:
\begin{description}
    \item[m1] $\Sigma_C \cap \varphi (\Sigma_C)= \{1_{G^0}\}$
    \item[m2] for all $\gamma \in G^{\circ}$ we have
    $\varphi(\gamma )\tau \gamma \notin \Sigma_C$.
\end{description}
Then there exists a free, mixed action of $G$ on $C \times C$,
hence $S= (C \times C)/G$ is a surface of general type isogenous
to a product of mixed type. More precisely, we have
\[ \gamma(x,y)=(\gamma x,\varphi (\gamma)y) \textit{ for } \gamma \in G^0
\]
\[ \tau'(x,y)=(y,\tau x).
\]
Conversely, every surface of general type isogenous to a product
of mixed type arises in this way.
\end{theo}
%----------------------------------------------------------------------%
Notice that $G^{\circ}$ is the subgroup of transformations not
exchanging the two factors of $C \times C$.
%-----------------------------------------------------------------------%
\begin{prop}\label{sottomix} Let $(C \times C)/G$ be a surface with $p_g=q=2$
isogenous to a product of mixed type. Then $B=C/G^{\circ}$ is a
curve of genus $g(B)=2$.
\end{prop}
%-------------------------------------------------------------------------%
\begin{proof}
From Proposition 3.15 \cite{cat00} we have:
\[ H^0(\Omega^1_S)=(H^0(\Omega^1_C) \oplus H^0(\Omega^1_C))^G=(H^0(\Omega^1_C)^{G^{\circ}} \oplus
H^0(\Omega^1_C)^{G^{\circ}})^{G/G^{\circ}}=
\]
\[=(H^0(\Omega^1_B) \oplus
H^0(\Omega^1_B))^{G/G^{\circ}}.
\]
Since $S$ is of mixed type, the quotient
$G/G^{\circ}=\mathbb{Z}_2$ exchanges the last summands, hence
$h^0(\Omega^1_S)=h^0(\Omega^1_B)=2$.
\end{proof}
%-----------------------------------------------------------------%
\begin{theo}\label{mainmix}
Let $S$ be a complex surface with $p_g=q=2$ and isogenous to a
product of mixed type. Then $S$ is minimal of general type and the
only possibility for the genus of the curve $C$ and the group $G$ is
given by the entry Mix of Table 1, moreover $G^{\circ}$ acts
freely on $C$.
\end{theo}
%-----------------------------------------------------------------%
\begin{proof}
From the fact that $| G | = (g(C)-1)^2$ and that $G^{\circ}$
is a subgroup of index $2$ in $G$ we have:
\begin{equation}\label{gzero}
| G^{\circ} | = \frac{(g(C)-1)^2}{2},
\end{equation}
hence $g(C)$ must be odd.
\\ Since $(C \times C)/ G^{\circ}$ is isogenous to a product of unmixed type and by Proposition \ref{sottomix}  we
have $g(C/G^{\circ})=2$, the Riemann-Hurwitz formula yields:
\[ 2g(C)-2=| G^{\circ} | \big(2+\sum^r_{i=1}(1-\frac{1}{m_i})\big),
\]
hence:
\[ 4=(g(C)-1)\big(2+\sum^r_{i=1}(1-\frac{1}{m_i})\big) \Rightarrow g(C)\leq 3 \Rightarrow g(C)=3
\mbox{ and } \sum^r_{i=1}(1-\frac{1}{m_i})=0.
\]
Then $ | G^0 | =2$ means $G^0=\mathbb{Z}_2$, $|G|=4$ and since (\ref{split})
is non-split $G=\mathbb{Z}_4$.
\end{proof}
%Moreover we see that the generating vector of $G^{\circ}$ in this
%case is unique, indeed if $\ZZ/2\ZZ=<x \mid x^2=1>$ the vector is:
%\[ a_1=x, \ \  b_1=1, \ \  a_2=1, \ \  b_2=1.
%\]
%We notice that the first case in Table 1 labelled with $UnMix$ was already given in \cite{zucc}.
%-------------------------------------------------------------------------------%
%-------------------------------------------------------------------------------%
\section{Isotrivial fibrations}\label{standiso}
Up to now we have considered only cases where a finite group $G$ acts freely
on a product of two curves $C \times F$, hence the quotient $(C \times F)/G$ is
smooth. Now we consider cases where $G$ does not act freely on $C \times F$.
Thus the quotient $(C \times F)/G$ is a normal surface, and we study its  desingularization.

\begin{defin} Assume that $T = ( C \times F )/G$, where $G$ is a finite group
of automorphisms of each factor $C$ and $F$, and acts diagonally on $C
\times F$. Consider the minimal resolution
$S$ of the singularities of $T$.
The holomorphic map $ f_C  \colon S \rightarrow C' : = C / G$ is called  a {\em
standard isotrivial fibration} if it is a relatively minimal fibration.
\end{defin}

\begin{rem} Let $S$ be a minimal surface of general type with $p_g(S)=q(S)=2$ and $S \rightarrow B$ be an isotrivial fibration with general fibre $F$, and $g(B)=1$, then by \cite[\S 1.1]{serr} $S$ is birational to $(C \times F)/G$, $B \cong C/G$ and by \eqref{eq.irregularity}
we have $g(C/G)=g(F/G)=1$. Consider the minimal desingularization $\sigma\colon S' \rightarrow (C \times F)/G$, then the holomorphic map $f_C  \colon S' \rightarrow C / G$ is a standard isotrivial fibration. Indeed suppose that there is a $-1$-curve $E$ in a fibre of $f_C$, then $\sigma(E)$ is a $-1$-curve in $(C \times F)/G$, but $(C \times F)/G \rightarrow C/G \times F/G$ is a finite map and $C/G \times F/G$ is a product of two elliptic curve, and this gives a contradiction. Thus $S \rightarrow B$ is birational to a standard isotrivial fibration, and we shall deal from now on only with standard isotrivial fibrations. 
\end{rem}

By abuse of notation we shall also denote by  $S \rightarrow T:= (C
\times F)/G$ a standard isotrivial fibration, and we shall refer to $S$ as a standard isotrivial fibration.

Let $S \rightarrow T:= (C
\times F)/G$ be a standard isotrivial fibration of general type, which is not  isogenous to a product of curves. To study the types of singularities of $T$, one looks first at the fixed points of the action of $G$ on each curve and at the stabilizers $H \subset G$ of each point on each curve. For this part we shall mainly follow \cite{mistretun}.

Let $C$ be a compact Riemann surface of genus $g \geq 2$ and let
$G \leq \Aut(C)$. For any $c \in G$ set $H:=\langle c \rangle$ and
define the set of fixed points by $c$ as:
\[ Fix_C(c)=Fix_C(H):=\{x \in C \mid cx=x \}.
\]

Let us look more closely to the action of an automorphism in a neighborhood of a fixed point. Let $\mathcal{D}$ be the unit disk and $c \in
\Aut(C)$ of order $m>1$ such that $cx=x$ for $x \in C$. Then there
is a unique complex primitive $m$-th root of unity $\xi$ such that any lift
of $c$ to $\mathcal{D}$ that fixes a point in $\mathcal{D}$ is
conjugate to the transformation $z \rightarrow \xi z$ in
$\Aut(\mathcal{D})$. We write $\xi_x(c)=\xi$ and we call $\xi^{-1}$
the \emph{rotation constant} of $c$ in $x$. Then for each integer
$q \leq m-1$ such that $(m,q)=1$ we define:
\[ Fix_{C,q}(c):=\{x \in Fix_C(c) \mid \xi_x(c)=\xi^q \},
\]
that is the set of fixed points of $c$ with rotation constant
$\xi^{-q}$. We have:
\begin{equation} Fix_C(c)=\biguplus_{\begin{array}{c} \scriptstyle q \leq m-1 \\
\scriptstyle (q,m)=1 \end{array}} Fix_{C,q}(c).
\end{equation}
\begin{lem}\cite[Lemma 10.4, Lemma 11.5]{breuer} Assume that we are in the situation of the Riemann
Existence Theorem \ref{RiemEx}, thus let $\mathcal{V}=(a_1,b_1, \dots
,a_{g'},b_{g'},c_1,\dots,c_{r})$ a generating vector of a
finite group $G$ of type $(g' \mid m_1, \dots, m_r)$. Let $c \in G
\setminus \{1\}$ be of order $m$, $H=\langle c \rangle$ and
$(q,m)=1$. Then:
\begin{equation}\label{eq.num.fix.pts}
|Fix_C(c)|=|N_G(H)| \sum_{\begin{array}{c} \scriptstyle 1 \leq i \leq r \\
\scriptstyle m \mid m_i \\
\scriptstyle H \sim_{G}\langle c^{m_i/m}_i\rangle \end{array}}
\frac{1}{m_i},
\end{equation}
and
\[ |Fix_{C,q}(c)|=|C_G(c)| \sum_{\begin{array}{c} \scriptstyle 1 \leq i \leq r \\
\scriptstyle m \mid m_i \\
\scriptstyle c \sim_{G}c^{qm_i/m}_i \end{array}}\frac{1}{m_i}.
\]
%where $N_G(x)$ is the normalizer of $x$ and $C_G(x)$ the
%centralizer of $x$.
\end{lem}
We need the following two Corollaries.
\begin{cor} Assume $c \sim_G c^q$. Then
$|Fix_{C,1}(c)|=|Fix_{C,q}(c)|$.
\end{cor}
\begin{cor}\label{cor.fix.centro} Let $c \in G$ with $|c|=2$ and $c \in Z(G)$, then:
\[ |Fix_C(c)|=|G| \sum_{\{i \mid c \in \langle c_i \rangle\}} \frac{1}{m_i}.
\]
\end{cor}

Let $S \rightarrow T:= (C
\times F)/G$ be a standard isotrivial fibration of general type as explained in \cite{serr}
paragraph (2.02) the stabilizer $H \subset G$ of a point $x \in
F$ is a cyclic group (see e.g.,
\cite{kar} Chap. III. 7.7), so, since the tangent representation is faithful on both factor, the only singularities that can occur
on $T$ are cyclic quotient singularities. More precisely, if $H$ is the
stabilizer of $x \in F$, then we have two cases.
If $H$ acts freely on $C$ then $T$ is
smooth along the scheme-theoretic fibre of $\sigma\colon T \rightarrow
F/G$ over $\overline{x} \in F/G$, and this fibre consists of the
curve $C/H$ counted with multiplicity $|H|$. Thus the
smooth fibres of $\sigma$ are all isomorphic to $C$. On the other
hand if a non-trivial element of $H$ fixes a point $y \in C$, then
$T$ has a cyclic quotient singularity in the point $\overline{(y,x)} \in
(C \times F)/G$. 
%---------------------------------------------------------------------------------%
 Let us briefly recall the definition of a cyclic quotient singularity.
\begin{defin}
Let $n$ and $q$ be natural numbers with $0 <q<n$ and $(n,q)=1$, and
let $\xi_n$ be a primitive $n-$th root of unity. Let the
action of the cyclic group $\mathbb{Z}_n = \left\langle \xi_n
\right\rangle$ on $\mathbb{C}^2$ be defined by $\xi_n \cdot
(x,y)=(\xi_nx,\xi^q_ny)$. Then we say that the analytic space
$\mathbb{C}^2/(\mathbb{Z}_n)$ has a \emph{cyclic quotient singularity}\index{Cyclic Singularity} of type
$\frac{1}{n}(1,q)$.
\end{defin}
We are interested in desingularization $\sigma\colon S
\rightarrow T$ of cyclic quotient singularities. The exceptional divisor $E$ on the minimal resolution of such a singularity is given by a
Hirzebruch-Jung string (see e.g., \cite{reid}, or \cite{BHPV}).
\begin{defin}
A Hirzebruch-Jung string is a union $E:=\cup^k_i E_i$ of smooth
rational curves $E_i$ such that:
\begin{itemize}
\item $E^2_i = -b_i  \leq -2 \textrm{ for all } i$,
\item $E_iE_j=1$  if
$| i-j | =1$,
\item $E_iE_j=0$ if  $| i-j | \geq 2$,
\end{itemize}
where the $b_i$'s are given by the
\emph{continued fraction} associated to $\frac{1}{n}(1,q)$. Indeed by the formula:
\[ \frac{n}{q}=b_1-\frac{1}{b_2-\frac{1}{\dots-\frac{1}{b_k}}}.
\]
\end{defin}
%Then the dual graph of $E$ is: {\setlength{\unitlength}{1.1cm}
%\begin{center}
%\begin{picture}(1,0.5)
%\put(0,0){\circle*{0.2}}
%\put(1,0){\circle*{0.2}} \put(0,0){\line(1,0){1}}
%\put(-0.3,0.2){\scriptsize $-b_1$}
%\put(0.70,0.2){\scriptsize $-b_2$}
%\put(2,0){\circle*{0.2}}
%\put(1,0){\line(1,0){0.2}}
%\put(1.3,0){\line(1,0){0.15}}
%\put(1.55,0){\line(1,0){0.15}}
%\put(1.8,0){\line(1,0){0.2}}
%\put(3,0){\circle*{0.2}} \put(2,0){\line(1,0){1}}
%\put(1.70,0.2){\scriptsize $-b_{k-1}$}
%\put(2.70,0.2){\scriptsize $-b_k$}
%\end{picture}         \hspace{2.5cm}
%\end{center}
%}

%\vspace{.5cm}
By abuse of notation we shall refer to $[b_1,\dots,b_k]$ as the continued fraction associated to $\frac{1}{n}(1,q)$.

These observations lead to the following theorem.
%-----------------------------------------------------------------------------%
\begin{theo}\cite[Theorem 2.1]{serr}\label{serano}Let $\sigma\colon S \rightarrow T:=(C \times F) / G$
be a standard isotrivial fibration and let us consider the natural
projection $\beta\colon S \rightarrow F/G$. Take any point over
$\overline{y} \in F/G$ and let $\Lambda$ denote the fibre of
$\beta$ over $\overline{y}$. Then:
\begin{enumerate}

    \item The reduced structure of $\Lambda$ is the union of an
    irreducible curve $Y$, called the central component of
    $\Lambda$, and either none or at least two mutually disjoint
    Hirzebruch-Jung strings, each meeting $Y$ at one point. These
    strings are in one-to-one correspondence with the branch
    points of $C \rightarrow C/H$, where $H$ is the stabilizer of
    $\overline{y}$.
    \item The intersection of a string with $Y$ is
    transversal, and it takes place at only one of the end
    components of the string.
    \item $Y$ is isomorphic to $C/H$, and has multiplicity
    $ |H| $ in $\Lambda$.
\end{enumerate}
Evidently, a similar statement holds if we consider the natural projection $\alpha\colon S \rightarrow C/G$.
\end{theo}

%--------------------------------------------------------------------------------%

We shall now determine the numerical invariants of isotrivial fibrations as we have done for surfaces isogenous to a product, in general the invariants will also depend on the singularities. The following Lemma, which derives from section $6$ of
\cite{barlow}, explains how to calculate them (for a more exhaustive treatment
of cyclic singularities of isotrivial surfaces with $\chi(\mathcal{O}_S)=1$ see
also \cite{mistretun}).
\begin{lem}\label{numeri}
Let $\sigma : S \rightarrow T:=C \times F / G$ be a standard isotrivial fibration. Then:
\begin{equation}\label{kiso}
 K^2_S=\frac{8(g(C)-1)(g(F)-1)}{| G |} +\sum_{x \in Sing (T)} h_x,
\end{equation}
\begin{equation}\label{eiso}
 e(S)=\frac{4(g(C)-1)(g(F)-1)}{| G |} +\sum_{x \in Sing (T)} e_x,
\end{equation}
where  $h_x$ and $e_x$ depend on the type of singularity in $x$.
If $x$ is a cyclic quotient singularity of type $\frac{1}{n}(1,q)$
then:
\[ h_x:= 2 - \frac{2+q+q'}{n}- \sum^k_{i=1}(b_i-2),
\]
\[ e_x:=k+1-\frac{1}{n},
\]
where $1 \leq q' \leq n-1$ and such that $qq' \equiv 1 mod \ n$, and
$b_i$ with ($1\leq i \leq k$) are the continued fractions data.
Moreover $e_x \geq \frac{3}{2}$.
\end{lem}
\begin{rem}\label{remuno} If $x \in T$ is a rational double point, i.e.,
$x$ is a singularity of type $\frac{1}{n}(1,n-1)$, we
have:
\[h_x = 0, \quad e_x = \frac{(n-1)(n+1)}{n}.
\]
In the case $x$ is a point of type $\frac{1}{3}(1,1)$, then the
continued fraction is given by $[3]$ hence we have:
\[h_x = -\frac{1}{3}, \quad e_x=\frac{5}{3}.
\]
In the case $x$ is a point of type $\frac{1}{4}(1,1)$, then the
continued fraction is given by $[4]$ hence we have:
\[h_x = -1, \quad e_x=\frac{7}{4}.
\]
In the case $x$ is a point of type $\frac{1}{5}(1,2)$, then the
continued fraction is given by $[3,2]$ hence we have:
\[h_x -\frac{2}{5}, \quad e_x=\frac{14}{5}.
\]
\end{rem}
\begin{rem}\label{remdue}Let us consider now a standard isotrivially fibred surface $S$ with $\chi(S)=1$, then by the Noether's formula:
\[ e(S)=12\chi(S)-K^2_S=12-K^2_S.
\]
Observe that together (\ref{kiso}) and (\ref{eiso}) yield:
\[ K^2_S=2e(S)-\sum_{x \in Sing (T)} (2e_x - h_x),
\]
Combining the above two formulas we get:
\begin{equation}\label{kappaiso}
K^2_S=8-\frac{1}{3}\sum_{x \in Sing (T)} (2e_x -h_x).
\end{equation}
Set $B_x:=2e_x -h_x$. Let us suppose now that $S$ is irregular,
hence by Debarre's inequality we have $K^2_S \geq 2$ (if
$p_g=q=2$ we even have $K^2_S \geq 4$). Since $S$ is smooth $\sum_{x \in Sing (T)}B_x \equiv 0$ $mod$ $3$, and combining these two facts gives the following upper bound:
\[ \sum_{x \in Sing (T)}B_x \leq 18.
\]
By the Bogomolov-Miyaoka-Yau inequality $K^2_S \leq 9$, this gives the
lower bound
\[ \sum_{x \in Sing (T)}B_x \geq -3.
\]
Now if $\sum_{x \in Sing (T)}B_x =0$ then $T$ is isomorphic to $S$, hence it is
smooth and $K^2_S=8$ and $S$ is isogenous to a product. If
$\sum_{x \in Sing (T)}B_x=-3$ then $S$ is a ball quotient by the
Miyaoka-Yau Theorem,
which is absurd. Hence $\sum_{x \in Sing (T)}B_x \geq 0$ and
\begin{equation}\label{kappotto}
K^2_S \leq 8.
\end{equation}
\end{rem}
We find a
basket of possible singularities for $T$ depending on $K^2_S$.
\begin{prop}\label{prop.list.sing} Let $\sigma\colon S \rightarrow T=(C \times F)/G$ be a standard isotrivial fibration with $\chi(S)=1$ and $p_g=q=2$. Then the possible singularities of $T$ are given in the following list.
\begin{itemize}
\item $K^2_S=6$:
\begin{enumerate}\renewcommand{\theenumi}{\it \roman{enumi}}
\item $2 \times \frac{1}{2}(1,1)$.
\end{enumerate}
\item $K^2_S=5$:
\begin{enumerate}\renewcommand{\theenumi}{\it \roman{enumi}}
\item $\frac{1}{3}(1,1)+\frac{1}{3}(1,2)$,
\item $2 \times \frac{1}{4}(1,1)$,
\item $3 \times \frac{1}{2}(1,1)$.
\end{enumerate}
\item $K^2_S=4$:
\begin{enumerate}\renewcommand{\theenumi}{\it \roman{enumi}}
\item $\frac{1}{4}(1,1)+\frac{1}{4}(1,3)$,
\item $2 \times \frac{1}{5}(1,2)$,
\item $\frac{1}{2}(1,1)+2 \times \frac{1}{4}(1,1)$,
\item $4 \times \frac{1}{2}(1,1)$.
\end{enumerate}
\end{itemize}
\end{prop}
We observe that this table is just a part of a more complete table
given in \cite[Proposition 4.1]{mistretun}, where the authors give
a list of all possible singularities for irregular standard
isotrivial fibrations with $\chi(S)=1$.

\begin{rem}
Recall from Remark \ref{singcong} that the singularities of $T$
arise from the points in $C \times F$ with non-trivial stabilizer;
since the action of $G$ on $C \times F$ is the diagonal one, it
follows that $\mathcal{S}'=(\Sigma (\mathcal{V}_1) \cap
\Sigma(\mathcal{V}_2)) \setminus \{1\}$ is the set of all
non-trivial stabilizers for the action of $G$ on $C \times F$.
Suppose that every element of $\mathcal{S}'$ has order $2$, then
we have that the singularities of $T$ are nodes, whose number is
given by:
\begin{equation}\label{eq.num.nodes}
\sharp \textrm{ Nodes}(T) = \frac{2}{|G|}\sum_{c \in
\mathcal{S}'}|Fix_{C}(c)||Fix_{F}(c)|,
\end{equation}
see e.g., \cite[\S 5]{pol1}.
\end{rem}

\begin{lem}\label{lemma.un.ram.data} Let $S$ be as in Theorem  \ref{mainisot} suppose $|Sing(T)|=2$ or $3$ and $g(F)=2$, then the covering $C \rightarrow C/G$ has only one ramification point.
\end{lem}
\begin{proof} Let us suppose that $C \rightarrow C/G$ has $r \geq 1$ ramification
points. Let $i \in \{1,\dots ,r\}$ and $\{m_i\}^r_{i=1}$ be the
branching data. Since $|Sing(T)|=2$ or $3$ the corresponding
Hirzebruch-Jung strings must belong to the same fibre of $S
\rightarrow C/G$, because by theorem \ref{serano} each fibre must
contain either none or at least two strings. It follows that, for
all $i$ except one there is a subgroup $ H \leq G$ isomorphic to
$\mathbb{Z}_{m_i}$, which acts freely on $F$. Now since $g(F)=2$
and by the Riemann-Hurwitz formula for the covering $F \rightarrow
F/H$ we have:
\[ 1=g(F)-1=m_i(g(F/H)-1),
\]
hence all the $m_i$'s except at most one divide $1$, therefore there is only one $m_i$, and so only one ramification point.
\end{proof}
\begin{theo}\label{mainisot}
Let $T:=(C \times F )/ G$ be a singular complex
surface and $S$ its minimal resolution of singularities such that $S$ has $p_g=q=2$.
Then $S$ is minimal of
general type and the only possibilities for $K^2_S$, $g(C)$, $g(F)$, the groups $G$, the branching data $\mathbf{m}$,
the types, and the numbers of singularities of $T$ are given in
Table 2.
\end{theo}
%%%%%%%%%%%%%%%%%%%%%%%%%%%%%%%%%%%%%%%%%%%%%%%%%%%%%%%%%%%%%%%%%%%%%%%%%%555
%%%%%%%%%%%%%%%%%%%%%%%%%%%%%%%%%%%%%%%%%%%%%%%%%%%%%%%%%%%%%%%%%%%%%%%%%%%%5
%
%%%%%%%%%%%%%%%%%%%%%%%%%%%%%%%%%%%%%%%%%%%%%%%%%%%%%%%%%%%%%%%%%%%%%%%%%%555
%%%%%%%%%%%%%%%%%%%%%%%%%%%%%%%%%%%%%%%%%%%%%%%%%%%%%%%%%%%%%%%%%%%%%%%%%%%%5
\begin{proof}

\textbf{Step 1} $S$ is minimal.
\\ First recall from Theorem \ref{Zucco} and (\ref{eq.irregularity}) that $S$ is of Albanese general type
and $C/G$ and $F/G$ are elliptic curves. If $E$ is a $(-1)$-curve
on $S$ then the image of $E$ in $T$ is rational curve. But $T \rightarrow C/G \times F/G$ is a finite map and $C/G \times F/G$ is a product of two elliptic curve, 
and this gives a contradiction, hence $S$
is minimal.

\textbf{Step 2} $4 \leq K^2_S \leq 6$.
\\ Since $S$ is minimal and irregular, by Debarre's inequality we have $K^2_S \geq 2 p_g= 4$, and by (\ref{kappotto})
we have $K^2_S \leq 8$. We see from equation (\ref{kappaiso}) that
if $K^2_S=8$ then $T$ is nonsingular, and this case cannot occur.
If $K^2_S=7$, then by (\ref{kappaiso}) we must have $\sum_{x \in
Sing (T)}B_x=3$, which means that $T$
can have only one singularity of type $\frac{1}{2}(1,1)$, but this contradicts Serrano's Theorem
\ref{serano}. Hence $4 \leq K^2_S \leq 6$.

\textbf{Step 3} consists of checking, once $K^2_S$ is fixed, if there is a standard isotrivial
surface $S \rightarrow T=(C \times F)/G$ with $p_g=q=2$ such that
$T$ has the prescribed singularities given in Proposition
\ref{prop.list.sing}.

\textbf{Case} $\mathbf{K^2_S=6}$.

In this case we have only a pair of singularities of type $\frac{1}{2}(1,1)$, hence $K^2_S=\frac{8(g(C)-1)(g(F)-1)}{|G|}$, and by Riemann-Hurwitz formula ($g(C/G)=1$) we have:
\begin{equation}\label{eq.r.h.ell} g(C)-1=\frac{|G|}{2}\sum^r_{i=1}(1-\frac{1}{m_i}).
\end{equation}
Combining the two formulas we have:
\[ \frac{3}{2}=(g(F)-1)\sum^r_{i=1}(1-\frac{1}{m_i}).
\]
Suppose $r=1$. Then $\frac{3}{2} \leq g(F)-1 \leq 3$, hence $3
\leq g(F) \leq 4$. By symmetry we can suppose w.l.o.g. $g(C) \leq
g(F)$, moreover we shall always assume this from now on.
\\ If $g(F)=4$ and $g(C)=4$, then both have one branching point of order $2$, and $|G|=12$.
 As we have already seen $A_4$ is the only group of order $12$ that is $(1 \mid 2)-$generated. Choose, for example, as generating vectors for $G$ for both coverings:
\[ a_{1,1}=a_{2,1}=(123), \ \ b_{1,1}=b_{2,1}=(124), \ \ c_{1,1}=c_{2,1}=(12)(34).
\]
We have $\mathcal{S}'=\{(12)(34),(13)(24),(14)(23)\}$. For all $c
\in \mathcal{S}'$ we have:
\[ |Fix_C(c)|=2, \textrm{  } |Fix_F(c)|=2, \]
so by equation \eqref{eq.num.nodes} $T$ has  $\frac{2 \cdot 2
\cdot 3}{6}=2$ nodes. Hence there exists $S$ and this gives the
last case of the Table 2.
\\ If $g(F)=4$ and $g(C)=3$, then $|G|=8$ and the branching data for $F$ and $C$ are respectively $(4)$ and $(2)$. However the commutators of $D_4$ and $Q_8$ have order $2$, hence neither group is $(1 \mid 4)-$generated, and this case is excluded.
\\ If $g(F)=3$, then $g(C)=3$ and $|G|=\frac{8 \cdot 2 \cdot 2}{6}=\frac{16}{3}$ which is absurd, and this case is impossible.

Suppose $r \geq 2$. Then $g(F)-1 \leq \frac{3}{2}$ hence $g(F)=2$
and this is a contradiction to Lemma \ref{lemma.un.ram.data}.

\textbf{Case} $\mathbf{K^2_S=5}$.

We have several cases according to Proposition \ref{prop.list.sing}.

\textit{Case (i)}. In this case we have two singularities, one of
type $\frac{1}{3}(1,1)$ and one of type $\frac{1}{3}(1,2)$. By
Remark \ref{remuno} we have $\sum h_x=-\frac{1}{3}$, hence we have
$K^2_S= \frac{8(g(C)-1)(g(F)-1)}{|G|} -\frac{1}{3}$, combining
this formula with (\ref{eq.r.h.ell}) we have:
\[ \frac{4}{3}=(g(F)-1)\sum^r_{i=1}(1-\frac{1}{m_i}).
\]
\\ Suppose $r=1$, then $2 < g(F) \leq 3$.
\\ If $g(F)=3$ and $g(C)=3$, then $|G|=6$, and the branching data
for both coverings are $(3)$. As we have seen $S_3$ is a non-abelian
group which is $(1 \mid 3)-$generated. Choose, for example, the
following generating vector for $G$ for both coverings:
\[ a_{1,1}=a_{2,1}=(12), \ \ b_{1,1}=b_{2,1}=(13), \ \
c_{1,1}=c_{2,1}=(123).
\]
We have $\mathcal{S}'=\{(123),(132)\}$ and for all $c \in
\mathcal{S}'$
\[ |Fix_{C,1}(c)| = |Fix_{C,2}(c)|= 1,
\]
\[ |Fix_{F,1}(c)| = |Fix_{F,2}(c)|= 1.
\]
So $C \times F$ contains exactly four points with non-trivial
stabilizer and for each of them the $G$-orbit has cardinality
$|G|/|\langle (123) \rangle|=2$. Hence $T$ contains precisely two
singular points and looking at the rotation constants we see that
it has the required singularities. Hence $S$ exists.
\\ If $r \geq 2$ then we have only the possibility $g(F)=2$ which is again a contradiction to Lemma \ref{lemma.un.ram.data}.

\textit{Case(ii)}. In this case we have two singularities of type
$\frac{1}{4}(1,1)$. By Remark \ref{remuno} $h_x=-1$,
which yields $K^2_S= \frac{8(g(C)-1)(g(F)-1)}{|G|} -2$, and combining
this formula with (\ref{eq.r.h.ell}) we have:
\[  \frac{7}{4}=(g(F)-1)\sum^r_{i=1}(1-\frac{1}{m_i}).
\]
Suppose $r=1$, then $ 3 \leq g(F) \leq 4$.
\\ If $g(F)=4$ and $g(C)=4$, then $|G|=\frac{8 \cdot 3 \cdot 3}{7}$, impossible.
\\ If $g(F)=4$ then $g(C)=3$ and $|G|=\frac{8\cdot 3 \cdot 2}{7}$, impossible.
\\ If $g(F)=3$ and $g(C)=3$, then $|G|=\frac{8 \cdot 2 \cdot 2}{7}$, impossible.
\\ Suppose $r \geq 2$ then the only possibility is $g(F)=2$ which is again a contradiction to Lemma \ref{lemma.un.ram.data}.

\textit{Case(iii)}. In this case we have three singularities of
type $\frac{1}{2}(1,1)$. By Remark \ref{remuno} we have $\sum
h_x=0$, which yields $K^2_S= \frac{8(g(C)-1)(g(F)-1)}{|G|}$, and
combining this formula with (\ref{eq.r.h.ell}) we obtain:
\[  \frac{5}{4}=(g(F)-1)\sum^r_{i=1}(1-\frac{1}{m_i}).
\]
Suppose $r=1$, then $ 2 < g(F) \leq 3$.
\\ If $g(F)=3$ and $g(C)=3$, then $|G|=\frac{8 \cdot 2 \cdot 2}{5}$, impossible.
\\ Suppose $r \geq 2$ then the only possibility is $g(F)=2$ which is again a contradiction to Lemma \ref{lemma.un.ram.data}.

\textbf{Case} $\mathbf{K^2_S=4}$.

We have several cases according to Proposition \ref{prop.list.sing}.

\textit{Case (i)}. In this case we have two singularities one
of type $\frac{1}{4}(1,1)$ and one of type $\frac{1}{4}(1,3)$. By
Remark \ref{remuno} we have $\sum h_x=-1$, hence we have $K^2_S=
\frac{8(g(C)-1)(g(F)-1)}{|G|} -1$. Combining this formula with
(\ref{eq.r.h.ell}) we obtain:
\[\frac{5}{4}=(g(F)-1)\sum^r_{i=1}(1-\frac{1}{m_i}).
\]
\\ Suppose $r=1$, then $2 < g(F) \leq 3$.
\\ If $g(F)=3$ and $g(C)=3$, then $|G|=\frac{8 \cdot 2 \cdot 2}{5}$, impossible.
\\ If $r \geq 2$ then we have only the possibility $g(F)=2$ which is again a contradiction to Lemma \ref{lemma.un.ram.data}.

\textit{Case(ii)}. In this case we have two singularities of type
$\frac{1}{5}(1,2)$. By Remark \ref{remuno} we have $\sum
h_x=-\frac{4}{5}$, which yields $K^2_S=
\frac{8(g(C)-1)(g(F)-1)}{|G|} -\frac{4}{5}$, combining this
formula with (\ref{eq.r.h.ell}) we have:
\[  \frac{6}{5}=(g(F)-1)\sum^r_{i=1}(1-\frac{1}{m_i}).
\]
Suppose $r=1$, then $ 2 < g(F) \leq 3$.
\\ If $g(F)=3$ and $g(C)=3$, then $|G|=\frac{8 \cdot 2 \cdot 2 \cdot 5}{24}$, impossible.
\\ Suppose $r \geq 2$, then the only possibility is $g(F)=2$ which is again a contradiction to Lemma \ref{lemma.un.ram.data}.

\textit{Case(iii)}. In this case we have three singularities one
of type $\frac{1}{2}(1,1)$ and two of type $\frac{1}{4}(1,1)$. By
Remark \ref{remuno} we have $\sum h_x=-2$, which yields $K^2_S=
\frac{8(g(C)-1)(g(F)-1)}{|G|}-2$, combining this formula with
(\ref{eq.r.h.ell}) we have:
\[ \frac{6}{4} =(g(F)-1)\sum^r_{i=1}(1-\frac{1}{m_i}).
\]
Suppose $r=1$, then $ 3 \leq g(F) \leq 4$.
\\ If $g(F)=4$ and $g(C)=4$, then $|G|=12$, and the branching data are $(2)$ for both coverings, but this contradicts the fact that we have singularities of type $\frac{1}{4}(1,1)$, hence this case is impossible.
\\  If $g(F)=4$ and $g(C)=3$, then $|G|=8$ and the branching data for $F$ and $C$ are respectively $(4)$, $(2)$. We have already seen that there is no non-abelian group of order $8$ which is $(1 \mid 4)-$generated, hence the case is excluded.
\\ If $g(F)=3$, then $g(C)=3$ and $|G|=\frac{8 \cdot 2 \cdot 2}{6}$, impossible.
\\ Suppose $r \geq 2$ then the only possibility is $g(F)=2$ which is again a contradiction to Lemma \ref{lemma.un.ram.data}.

\textit{Case(iv)}. In this case we have four singularities of type $\frac{1}{2}(1,1)$. By Remark \ref{remuno} we have $\sum h_x=0$, which yields $K^2_S= \frac{8(g(C)-1)(g(F)-1)}{|G|}$, and combining this formula with (\ref{eq.r.h.ell}) we have:
\[  1=(g(F)-1)\sum^r_{i=1}(1-\frac{1}{m_i}).
\]
Suppose $r=1$, then $ 2 < g(F) \leq 3$.
\\ If $g(F)=3$ and $g(C)=3$, then $|G|=8$, and the branching data are $(2)$ for both coverings. The groups $D_4$ and $Q_8$ are $(1 \mid 2)-$generated, choose, for example, as generating vector for $Q_8$ for both coverings:
\[ a_{1,1}=a_{2,1}=i, \ \ b_{1,1}=b_{2,1}=j, \ \ c_{1,1}=c_{2,1}=-1,
\]
and for $D_4$ for both coverings:
\[ a_{1,1}=a_{2,1}=x, \ \ b_{1,1}=b_{2,1}=y, \ \
c_{1,1}=c_{2,1}=x^2.
\]
We have in both cases $\mathcal{S}'=Z(G)\setminus \{1\}$, and by
Corollary \ref{cor.fix.centro} for $c \in \mathcal{S}'$
\[ |Fix_{C}(c)| = |Fix_{F}(c)|= 4.
\]
Then by equation \eqref{eq.num.nodes} $T$ has exactly $\frac{2
\cdot 4 \cdot 4}{8}=4$ nodes in both cases. Hence $S$ exists.
\\ Suppose $r \geq 2$ then the only possibility is $g(F)=2$ and $r=2$.
In this case there are more than three singularities so Lemma
\ref{lemma.un.ram.data} does not apply, hence $g(C)=2$, $|G| =2$
and both coverings have branching data $(2,2)$. Let $x$ be the
generator of $G$, we have $\mathcal{S}'=\{x\}$, and Corollary
\ref{cor.fix.centro} implies
\[ |Fix_{C}(x)| = |Fix_{F}(x)|= 2.
\] 
Then by
equation \eqref{eq.num.nodes} $T$ has exactly $\frac{2 \cdot 2
\cdot 2}{2}=4$ nodes. This yields the first case in the table.

\end{proof}
We notice that the first case in Table 2 was already given in \cite{zucc}.
%\begin{rem}
%In the end we can say, that our classification completes the one given in \cite{zucc} adding: the cases where the %groups acting on the product
%of curves are not abelian and the surfaces are of Albanese general type, and the mixed case.
%\end{rem}
%-------------------------------------------------------------------------------%
%-------------------------------------------------------------------------------%
\section{Moduli Spaces}\label{moduli}
By a celebrated Theorem of Gieseker (see \cite{gie}), once the two
invariants of a minimal surface $S$ of general type, $K^2_S$ and
$\chi(S)$, are fixed, then there exists a quasiprojective moduli
space $\mathcal{M}_{K^2_S,\chi(S)}$ of minimal smooth complex
surfaces of general type with those invariants, and this space
consists of a finite number of connected components. The union
$\mathcal{M}$ over all admissible pairs of invariants ($K^2,
\chi$) of these spaces is called the \emph{moduli space of
surfaces of general type}.

In \cite{cat00}, Catanese studied the moduli space of surfaces
isogenous to a product of curves (see Theorem 4.14). As a
result, one obtains that the moduli space of surfaces isogenous to
a product of curves $\mathcal{M}_{(G,(\tau_1,\tau_2))}$ with fixed invariants: a finite group
$G$ and types $(\tau_1,\tau_2)$  (where the types $\tau_i:=(g'_i |
m_{i,1},\dots,m_{i,r_i})$, for $i=1,2$, are defined in
\ref{genrvect}) for the unmixed case (while only
$G$ and one type $\tau$ in the mixed case), consists of a finite
number of irreducible connected components of $\mathcal{M}$.

The surfaces we are studying are quotients of products of curves and to study
their moduli space one has to look first at the moduli space of
Riemann surfaces.

Let $\mathcal{M}_{g,r}$ denote the \emph{moduli space of Riemann
surfaces of genus $g$ with $r$ ordered marked points}. By
permuting the marked points on the Riemann surfaces, the
permutation group $S_r$ acts naturally on this space; the moduli
space $\mathcal{M}_{g,[r]} = \mathcal{M}_{g,r}/S_r$ classifies the
Riemann surfaces of genus $g$ with $r$ unordered marked points. By
Teichm\"{u}ller theory these spaces are quotient of a contractible
spaces $\mathcal{T}_{g,r}$ of complex dimension $3g-3+r$, if $g'=0$ and $r \geq 3$, or $g'=1$ and $r\geq 1$ or $g'\geq 2$, called
the \emph{Teichm\"{u}ller spaces}, by the action of discrete
groups called the \emph{full mapping class groups} $Map_{g,[r]}$.

In \cite[Theorem 1.3]{ingridfab} is given a method to calculate the
number of connected components of the moduli spaces
$\mathcal{M}_{(G,(\tau_1,\tau_2))}$ of surfaces isogenous to a product of unmixed type using Teichm\"{u}ller theory, while
in \cite[Proposition 5.5]{infabfriz} is treated the mixed case.

Notice that the dimension of the space
$\mathcal{M}_{(G,(\tau_1,\tau_2))}$, with types $\tau_i:=(g'_i |
m_{i,1},\dots,m_{i,r_i})$ for $i=1,2$, is precisely
$dim\mathcal{M}_{(G,(\tau_1,\tau_2))}= 3g'_1-3+r_1+3g'_2-3+r_2$, while in the mixed case, if $\tau=(g'|m_1, \cdots ,m_r)$, then $dim
\mathcal{M}_{(G,\tau)}=3g'-3+r$. This is enough to determine the
numbers in the column $dim$ of Table 1.

\begin{defin} Let $M$ be a differentiable manifold, then the \emph{mapping
class group}\index{Mapping Class Group} (or Dehn group) of $M$ is the group:
\[ \Map(M):= \pi_0(\Diff^+(M))=\Diff^+(M)/\Diff^0(M),
\]
where $\Diff^+(M)$ is the group of orientation preserving
diffeomorphisms of $M$ and $\Diff^0(M)$ is the subgroup of
diffeomorphisms of $M$ isotopic to the identity.
\\ If $M$ is a compact complex curve of genus $g'$ we will use the following notations:
\begin{enumerate}
    \item We denote the mapping class group of $M$ by
    $\Map_{g'}$.
    \item If we consider $r$ points $p_1,\dots, p_r$ on $M$ we define:
\[ Map_{g',[r]}:=\pi_0(\Diff^+(M-\{p_1,\dots,p_r\})),
\]
and this is known as the \emph{full mapping class group}.
\end{enumerate}
\end{defin}

There is an easy way to present the mapping class group of a
curve using half twists and Dehn twists see e.g., \cite{catdd,deh}.
%\begin{prop} The mapping class group
%$Map_{0,[r]}$ is generated by $\sigma_1,\dots,\sigma_r $ with the
%following relations:
%\[
%\sigma_i\sigma_{i+1}\sigma_i=\sigma_{i+1}\sigma_{i}\sigma_{i+1},
%\]
%\[ \sigma_i\sigma_{j}=\sigma_j\sigma_{i}, \ \ \ \ \textrm{if} \ \ \ \ |
%i-j | \geq 2,
%\]
%\[
%\sigma_{r-1}\sigma_{r-2}\dots\sigma^2_1\dots\sigma_{r-2}\sigma_{r-1}=1.
%\]
%\end{prop}
%For a proof of the above Theorem see for example \cite[Theorem 1.11]{Bir}. 

%We want to give a similar presentation for a group $Map_{2}$. 
\begin{theo} The group $\Map_2$ is generated by the Dehn
twists with respect to the five curves in the figure:
\begin{center}
%\begin{minipage}[b]{8.2cm}
\includegraphics[width=8cm]{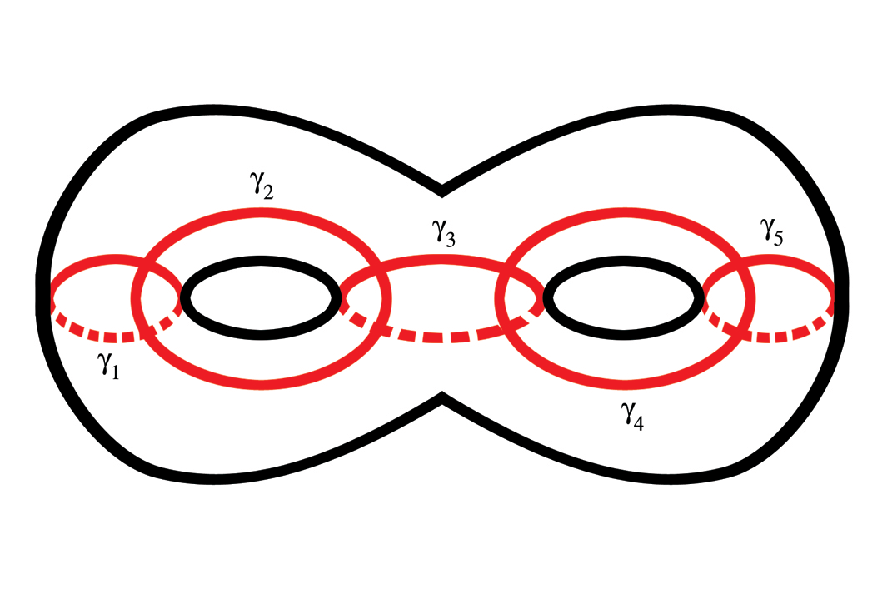}
%\end{minipage}
Figure 1.
\end{center}
The corresponding relations are the following:
\begin{enumerate}
    \item $t_{\gamma_i}t_{\gamma_j}=t_{\gamma_j}t_{\gamma_i}$ if $|i-j|\geq 2$,   $1 \leq i,j \leq 5$,
    \item $t_{\gamma_i} t_{\gamma_{i+1}} t_{\gamma_i} = t_{\gamma_{i+1}}t_{\gamma_{i}}t_{\gamma_{i+1}}$,  $1 \leq i \leq 4$,
    \item $(t_{\gamma_1}t_{\gamma_2}t_{\gamma_3}t_{\gamma_4}t_{\gamma_5})^6=1$,
    \item $(t_{\gamma_1}t_{\gamma_2}t_{\gamma_3}t_{\gamma_4}t_{\gamma_5}^2t_{\gamma_4}t_{\gamma_3}t_{\gamma_2}t_{\gamma_1})^2=1$,
    \item $[t_{\gamma_1}t_{\gamma_2}t_{\gamma_3}t_{\gamma_4}t_{\gamma_5}^2t_{\gamma_4}t_{\gamma_3}t_{\gamma_2}t_{\gamma_1}, t_{\gamma_i}]=1$, $1 \leq i \leq 5$.
\end{enumerate}
\end{theo}
A proof of the above theorem can be found in \cite[Theorem 4.8]{Bir}.

%-------------------------------------------------------------------------------%
%We have the following exact sequence of groups (cf. \cite{Macl} \S 4):
%\begin{equation}\label{ex.seq.mod.aut} 1 \longrightarrow \Inn(\Gamma) \longrightarrow \Aut^+(\Gamma) \longrightarrow \Out^+(\Gamma) \longrightarrow 1.
%\end{equation}
%%%%%%%%%%%%%%%%%%%%%%%%%%%%%%%%%%%%%%%%%%%%%%%%%%%%%%%%%%%%%%%%%%%%%%%%%%%%%%%%%%%%%%%%%%%%
%\begin{defin} Let $C$ be a complex Riemann surface of genus $g$, and $\Hom^{+}(C)$ be its group of orientation preserving homeomorphisms.
%\end{defin}
%Following \cite{Broug} section 2 we have

%%%%%%%%%%%%%%%%%%%%%%%%%%%%%%%%%%%%%%%%%%%%%%%%%%%%%%%%%%%%%%%%%%%%%%%%%%%%%%%%%%%%%%%%%%%%
%bla bla
%Let $C$ be a complex curve of genus $g'\geq 2$, consider the representation of $C$ as a quotient space $\mathcal{H}/\Pi_{g'}$ of the upper half plane $\mathcal{H}$ by the free action of a Fuchsian group $\Pi_g \cong  \pi_g(C)$. Then an  automorphism of $
%-------------------------------------------------------------------------------%
\begin{theo}\label{out}
Let $\Gamma=\Gamma(g' \mid m_1,\dots,m_r)$ be
an orbifold surface group and $Out^+(\Gamma)$ the group of orientation preseving outer automorphisms of $\Gamma$. Then there is an isomorphism of groups:
\[ \Out^+(\Gamma) \cong \Map_{g',[r]}.
\]
\end{theo}
This is a classical result cf. e.g., \cite[Theorem 2.2.1]{schn} and \cite[\S 4]{Macl}.

Moreover let $G$ be a finite group $(g'
\mid m_1,\dots,m_r)-$generated, then the action of $Out^+(\Gamma) \cong \Map_{g',[r]}$ on the generators of $\Gamma(g'
\mid m_1,\dots,m_r)$ induces an action on the generating vectors of $G$ via composition with admissible epimorphisms.

%following \cite{Broug} Section 2, there is a natural action
%of $\Aut^+(\Gamma)$ on the set of admissible epimorphisms from $\Gamma$ to $G$ given
%$Epi(\pi_1(C),\Gamma, G):=\{\theta: \Gamma \rightarrow G \textrm{ : }\theta \textrm{ surjective, } ker(\theta)$
%$\textrm{ torsion free}\}$ given
%by:
%\[ \eta(\theta)=\theta \circ \eta^{-1}.
%\]
%Notice that taking a lift in $\Aut^+(\Gamma)$ of an element in $\Out^+(\Gamma)$ we have an %action of  $\Out^+(\Gamma)$ on the set of admissible epimorphisms well defined up to inner %automorphisms. This action induces an action on the set of systems
%of generators of $G$.
%$\Out^+(\Gamma)$ on the set of admissible epimorphisms.

\begin{defin}\label{humov} Let $G$ be a finite group $(g'
\mid m_1,\dots,m_r)-$generated. 
%The group $\Map_{g',[r]}$ acts on the
%set of systems of generators of $G$, where the action is induced by the
%action on admissible epimorphisms. 
If two generating vectors
$\mathcal{V}_1$ and $\mathcal{V}_2$ are in the same $\Map_{g',[r]}$-orbit, we say that
they are related by a \emph{Hurwitz move}  (or are \emph{Hurwitz equivalent}).
\end{defin}
Now we calculate the Hurwitz moves on the generating vectors of type $(2 \mid -)$.
%-------------------------------------------------------------------------------%
%\begin{prop} Up to inner automorphism, the action of
%$\Map_{0,[r]}$ on $\Gamma (0 \mid m^r)$ is given by:
%\[
%\sigma_i: \left\{
%\begin{array}{rl} \gamma_i & \rightarrow \gamma_{i+1} \\
%\gamma_{i+1}  & \rightarrow \gamma^{-1}_{i+1}\gamma_i\gamma_{i+1} \\
%\gamma_j  & \rightarrow \gamma_{j} \textrm{ if } j \neq i, i+1.
%\end{array}
%\right.
%\]
%\end{prop}
%Cf. \cite{catdd} section 5.
%\begin{rem} Now we induce the action on the spherical generating vector.
%We see, that if we have different $m_i$, we have to consider
%permutations of the marked points in order not to mix the
%branching data, which correspond to the order of the elements. To
%solve this problem we consider Definition
%\ref{genrvect} with $(ii)$ \textbf{B}.
%\end{rem}
%According to the Remark above we have the following Corollary.
%-------------------------------------------------------------------------------%
%\begin{cor} Let $G$ be a finite group and let
%$\mathcal{V}=(c_1,\dots,c_r)$ be a spherical unordered generating vector for $G$ of type $\tau=(m_1,\dots,m_r)$. Then the Hurwitz
%moves on the set of spherical generating vectors of $G$ of type $\tau$ are generated by:
%\[
%\sigma_i: (c_1,\dots,c_r) \longmapsto (c'_1,\dots,c'_r),
%\]
%where
%\[
%\begin{array}{rl} c'_i & = c_{i+1}, \\
%c'_{i+1} & = c^{-1}_{i+1}c_ic_{i+1} ,\\
%c'_j  & = c_{j} \textrm{ if } j \neq i, i+1.
%\end{array}
%\]

%\end{cor}
%Let us look at the second case now.
%-------------------------------------------------------------------------------%
\begin{prop}\label{prop.hrwz.2} Up to inner automorphism, the action of
$\Map_{2}$ on $\Gamma (2 \mid -)$ is given by:
\[
t_{\gamma_2}: \left\{
\begin{array}{rl}
\alpha_1 & \rightarrow \alpha_1 \\
\beta_{1}  & \rightarrow \beta_{1}\alpha_{1} \\
\alpha_2 & \rightarrow \alpha_2 \\
\beta_2 & \rightarrow \beta_2 \\
\end{array}
\right. \ \ \ \
t_{\gamma_1}: \left\{
\begin{array}{rl}
\alpha_1 & \rightarrow \alpha_1\beta^{-1}_1 \\
\beta_{1}  & \rightarrow \beta_{1} \\
\alpha_2 & \rightarrow \alpha_2 \\
\beta_2 & \rightarrow \beta_2 \\
\end{array}
\right.
\]
\[
t_{\gamma_5}: \left\{
\begin{array}{rl}
\alpha_1 & \rightarrow \alpha_1 \\
\beta_{1}  & \rightarrow \beta_{1} \\
\alpha_2 & \rightarrow \alpha_2\beta^{-1}_2 \\
\beta_2 & \rightarrow \beta_2 \\
\end{array}
\right. \ \ \ \
t_{\gamma_4}: \left\{
\begin{array}{rl}
\alpha_1 & \rightarrow \alpha_1 \\
\beta_{1}  & \rightarrow \beta_{1} \\
\alpha_2 & \rightarrow \alpha_2 \\
\beta_2 & \rightarrow \beta_2\alpha_2 \\
\end{array}
\right.
\]
\[
t_{\gamma_3}: \left\{
\begin{array}{rl}
\alpha_1 & \rightarrow \alpha_1 x^{-1} \\
\beta_{1}  & \rightarrow x \beta_{1}x^{-1} \\
\alpha_2 & \rightarrow x\alpha_2 \\
\beta_2 & \rightarrow \beta_2. \\
\end{array}
\right.
\]
where $\alpha_1$, $\alpha_2$, $\beta_1$ and $\beta_2$ are the generators of
$\Gamma(2 \mid -)$ and $x=\beta^{-1}_2\alpha_1\beta_1\alpha^{-1}_1=\alpha_2\beta^{-1}_2\alpha^{-1}_2\beta_1$.
\end{prop}
\begin{proof}
One notices that a Riemann surface of genus $2$ is a
connected sum of two tori. Then one can use the results given in
\cite[Proposition 1.10]{pol1} to calculate the Dehn twists about
the curves $\gamma_1$, $\gamma_2$, $\gamma_4$, $\gamma_5$ of Figure 1,
considering the action on the two different tori. This gives the
actions $t_{\gamma_1}$, $t_{\gamma_2}$, $t_{\gamma_4}$ and
$t_{\gamma_{5}}$.
\\ Then the only Dehn twist left to calculate is the one with respect to the curve
$\gamma_3$ as in Figure 2.
%---------------------------------------------%
\begin{center}
%\begin{figure}
%\begin{minipage}[b]{8.2cm}
\includegraphics[width=8cm]{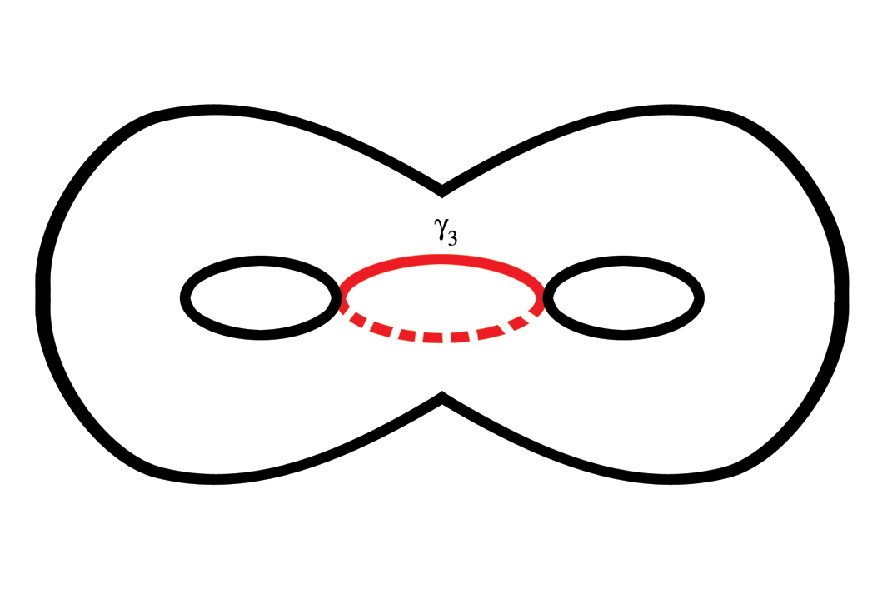}
%\end{minipage}
Figure 2.
%\end{figure}
\end{center}
%---------------------------------------------%
Choose the generators of the fundamental group as in Figure 3:
%---------------------------------------------%
\begin{center}
%\begin{minipage}[b]{8.2cm}
\includegraphics[width=8cm]{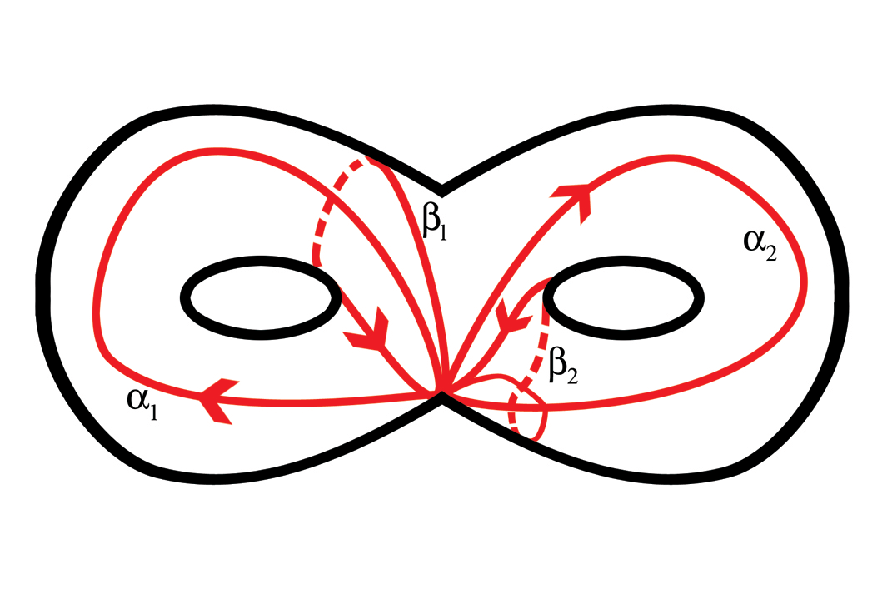}
%\end{minipage}
Figure 3.
\end{center}
%--------------------------------------------%
One sees that the only curves which have to be twisted are
$\alpha_1$, $\beta_1$ and $\alpha_2$ because the other is disjoint from $\gamma_3$
. In Figure 4 one sees the Dehn twist of $\alpha_1$ with respect to
$\gamma_3$. Following the curve one constructs the image of $\alpha_1$
under the map $t_{\gamma_3}$.
%--------------------------------------------%
\begin{center}
%\begin{minipage}[b]{8.2cm}
\includegraphics[width=8cm]{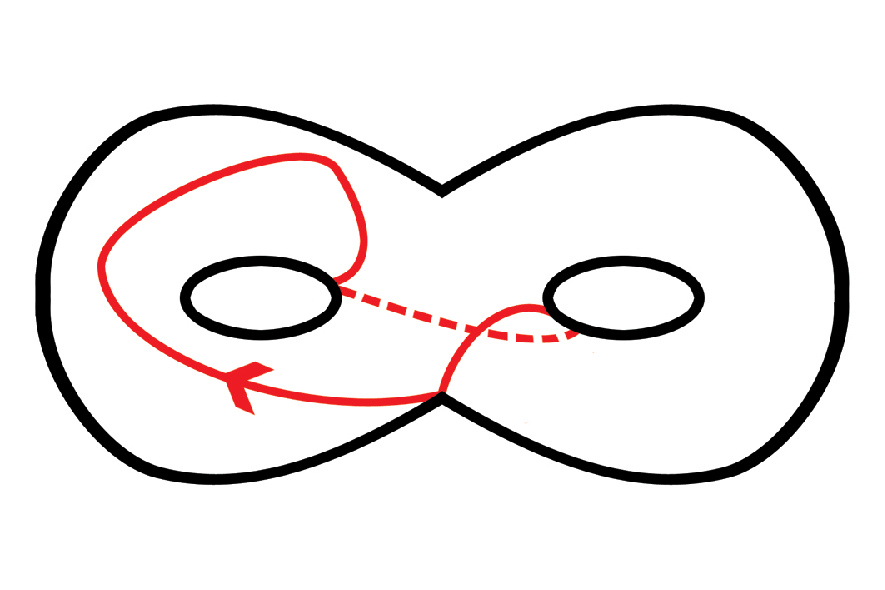}
%\end{minipage}
Figure 4.
\end{center}
%--------------------------------------------%
In Figure 5 we give the Dehn twist of $\beta_1$ with respect to
$\gamma_3$.
%--------------------------------------------%
\begin{center}
%\begin{minipage}[b]{8.2cm}
\includegraphics[width=8cm]{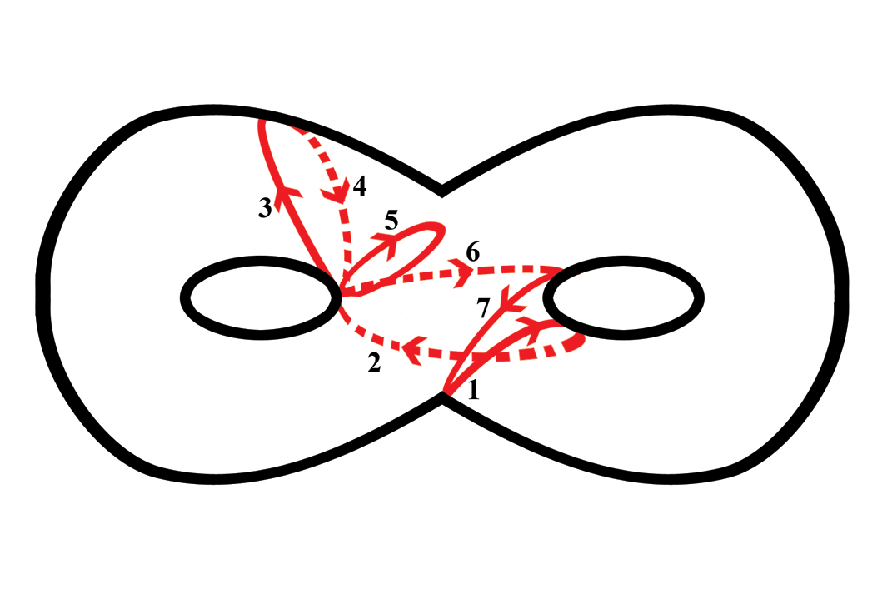}
%\end{minipage}
Figure 5.
\end{center}
%-----------------------------------------------%
In the last Figure we give the Dehn twist of $\alpha_2$ with respect to
$\gamma_3$ which completes the proof.
%-----------------------------------------------%
\begin{center}
%\begin{minipage}[b]{8.2cm}
\includegraphics[width=8cm]{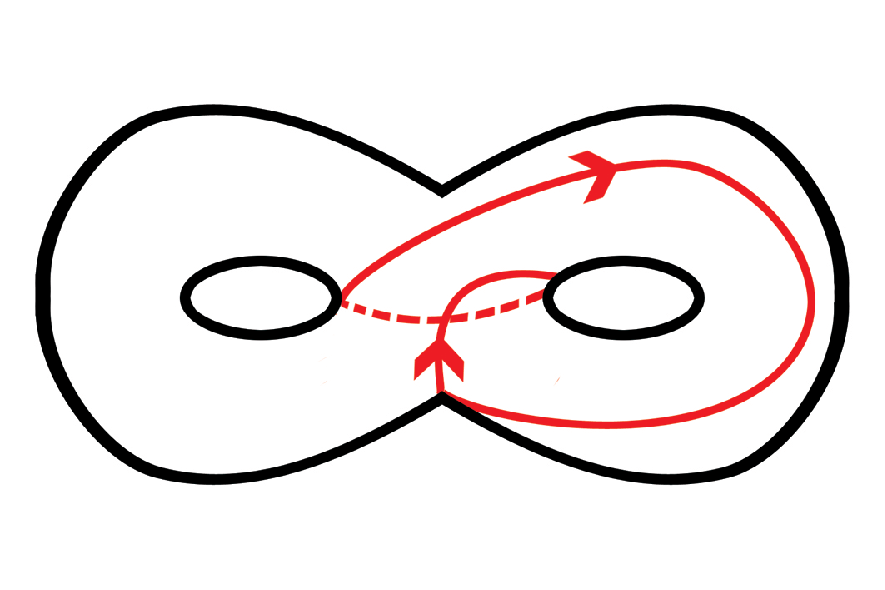}
%\end{minipage}
Figure 6.
\end{center}
%-----------------------------------------------%
\end{proof}
%-------------------------------------------------------------------------------%
\begin{cor} Let $G$ be a finite group and let
$\mathcal{V}=(a_1,b_1,a_2,b_2)$ be a generating vector for $G$ of type
$\tau=(2 \mid -)$. Then the Hurwitz moves on the set of generating vectors of $G$ of type $\tau$ are generated by:
\[
\mathbf{1}: \left\{
\begin{array}{rl}
a_1 & \rightarrow a_1 \\
b_{1}  & \rightarrow b_{1}a_{1} \\
a_2 & \rightarrow a_2 \\
b_2 & \rightarrow b_2 \\
\end{array}
\right. \ \ \ \ \mathbf{2}: \left\{
\begin{array}{rl}
a_1 & \rightarrow a_1b^{-1}_1 \\
b_{1}  & \rightarrow b_{1} \\
a_2 & \rightarrow a_2 \\
b_2 & \rightarrow b_2 \\
\end{array}
\right.
\]
\[
\mathbf{3}: \left\{
\begin{array}{rl}
a_1 & \rightarrow a_1 \\
b_{1}  & \rightarrow b_{1} \\
a_2 & \rightarrow a_2b^{-1}_2 \\
b_2 & \rightarrow b_2 \\
\end{array}
\right. \ \ \ \
\mathbf{4}: \left\{
\begin{array}{rl}
a_1 & \rightarrow a_1 \\
b_{1}  & \rightarrow b_{1} \\
a_2 & \rightarrow a_2 \\
b_2 & \rightarrow b_2a_2 \\
\end{array}
\right.
\]
\[
\mathbf{5}: \left\{
\begin{array}{rl}
a_1 & \rightarrow a_1 x^{-1} \\
b_{1}  & \rightarrow x b_{1}x^{-1} \\
a_2 & \rightarrow xa_2 \\
b_2 & \rightarrow b_2. \\
\end{array}
\right.
\]
where $x=b^{-1}_2a_1b_1a^{-1}_1=a_2b^{-1}_2a^{-1}_2b_1$.
\end{cor}
\begin{proof}
This follows directly from Proposition \ref{prop.hrwz.2}.
\end{proof}
%-------------------------------------------------------------------------------%
For the Hurwitz moves for types $(0 \mid m_1, \cdots ,m_r)$, $(1 \mid 1)$ and $(1 \mid 2)$ we refer to \cite{pol1} for proofs.
If $\mathcal{V}:=\{c_1, \ldots, c_r\}$ is an unordered generating vector of type $\tau:=(0 \; | \; m_1,
\ldots, m_r)$ then the Hurwitz moves on the set of unordered generating vectors of $G$ of type $\tau$ are generated, for $1 \leq \mathbf{i} \leq r-1$, by
\begin{equation*}
\mathbf{i}: \left \{
\begin{array}{ll}
c_i & \lr c_{i+1} \\
c_{i+1} & \lr c_{i+1}^{-1} c_i c_{i+1} \\
c_j & \lr c_j  \quad \textrm{if} \; j \neq i,\; i+1.
\end{array} \right.
\end{equation*}
If $\mathcal{V}:=\{a_1, b_1, c_1\}$ is of type $(1 \; | \; m)$
then the Hurwitz moves are generated by
\begin{equation*}
\mathbf{1} \colon \left\{ \begin{array}{ll}
c_1 \lr c_1 \\
a_1 \lr a_1 \\
b_1 \lr b_1a_1 \\
\end{array} \right.
\quad \mathbf{2} \colon \left\{ \begin{array}{ll}
c_1 \lr c_1 \\
a_1 \lr a_1b_1^{-1} \\
b_1 \lr b_1. \\
\end{array} \right.
\end{equation*}
If $\mathcal{V}:=\{a_1, b_1, \, c_1, c_2\}$ is of type $(1 \;
| \; m^2)$
 then the Hurwitz moves are generated by
\begin{equation*}
\begin{split}
& \mathbf{1} \colon \left\{ \begin{array}{ll}
  c_1 \lr c_1 \\
  c_2 \lr c_2 \\
  a_1 \lr a_1 \\
  b_1 \lr b_1a_1
\end{array} \right.
\quad  \quad \quad \quad \quad \quad \quad \mathbf{2} \colon \left\{
\begin{array}{ll}
  c_1 \lr c_1 \\
  c_2 \lr c_2\\
  a_1 \lr a_1b_1^{-1} \\
  b_1 \lr b_1
\end{array} \right.
\\
& \mathbf{3} \colon \left\{ \begin{array}{ll}
  c_1 \lr c_1 \\
  c_2 \lr a_1b_1^{-1}a_1^{-1}c_2a_1b_1a_1^{-1}\\
  a_1 \lr b_1^{-1}c_1a_1 \\
  b_1 \lr b_1
\end{array} \right.
\quad  \; \mathbf{4} \colon \left\{ \begin{array}{ll}
  c_1 \lr b_1^{-1}a_1^{-1}c_2 a_1b_1 \\
  c_2 \lr a_1^{-1}b_1^{-1}c_1 b_1a_1 \\
  a_1 \lr a_1^{-1} \\
  b_1 \lr b_1^{-1}.
\end{array} \right.
\end{split}
\end{equation*}
%------------------------------------------------------------------------------%
The following
theorem is a natural generalization of \cite[Theorem 1.3]{ingridfab}.
%-------------------------------------------------------------------------------%
\begin{theo}\label{Fabmain} Let $S$ be a surface isogenous to a
product of unmixed type with $p_g=q=2$. Then to $S$ we attach its
finite group $G$ (up to isomorphism) and the equivalent class of
a pair of disjoint unordered generating vectors
$(\mathcal{V}_1,\mathcal{V}_2)$ of type $((2 \mid -),(0\mid
\mathbf{m}))$ (or  $((1 \mid \mathbf{n}_1),(1\mid \mathbf{n}_2))$) of
$G$, under the equivalent relation generated by:
\begin{enumerate}
    \item Hurwitz equivalence and $\Inn(G)$ on $\mathcal{V}_1$,
    \item Hurwitz equivalence and $\Inn(G)$ on $\mathcal{V}_2$,
    \item simultaneous conjugation of $\mathcal{V}_1$ and
    $\mathcal{V}_2$ by an element $\lambda \in Aut(G)$, i.e., we let
    $(\mathcal{V}_1,\mathcal{V}_2)$ be equivalent to
    $(\lambda(\mathcal{V}_1),\lambda(\mathcal{V}_2))$.
\end{enumerate}
Then two surfaces $S$ and $S'$ are deformation equivalent if and
only if the corresponding pairs of generating vectors are in the
same equivalence class.
\end{theo}
\begin{proof}
We can use the same argument of \cite[Theorem 1.3]{ingridfab},
one has only to substitute the braid actions with the action of the
appropriate mapping class group.
\end{proof}
%--------------------------------------------------------------------------------%

To calculate the number of connected components is then the same
as to calculate the number of all possible generating pairs modulo the
equivalence relation defined above. Since this task may be too
hard to be achieved by hand with S. Rollenske we developed a
program in GAP4 which calculates the number of these pairs. The appendix by
Rollenske is devoted to explain how the program works, for the script of the program see \cite{io}. The result
of this computation is the last column of Table 1.

\begin{rem}\label{rem.almost.right} Notice that in the program is not implemented the action generated by the mapping class group and the group of inner automorphisms of the group $G$, but only the action of the mapping class group, hence we do not act with full group prescribed by Theorem \ref{Fabmain}. 

However this does not effect the result above. Indeed in the cases where there is only one orbit this is not a problem. It is not a problem neither if there are two orbits and the group $G$ is abelian, because in this case the inner automorphisms act trivially. 

But there are three cases where the group $G$ is not abelian and we have two orbits. In all three cases the two pairs of generating vectors are of the form $(\mathcal{V}_1,\mathcal{V}_2)$ and $(\mathcal{V}_1,\mathcal{W}_2)$. Then we used the program to calculate the orbits, only on the right side of the pairs, under the action of the group generated by the mapping class group and the group of automorphisms of $G$. Notice that this group contains the group generated by the mapping class group and the group of inner automorphisms. As result we have two orbits, hence we have two orbits also for the action of the group prescribed in Theorem \ref{Fabmain}. 

Moreover in \cite{io} are exhibited the pairs $(\mathcal{V}_1,\mathcal{V}_2)$ of generating vectors which give the surfaces
isogenous to a product of curves of unmixed type with
$p_g=q=2$ given in Table 1.
\end{rem}

For the mixed case we notice that there is only one connected
component of dimension $3$ of the moduli space corresponding to
the item labelled by Mix in Table 1. This comes directly from
the proof of Theorem \ref{mainmix} and from \cite[Proposition
5.5]{infabfriz} adapted to this case. Indeed let us denote by
$\mathcal{M}_{(G,\tau)}$ the moduli space of isomorphism classes
of surfaces isogenous to a product of curves of mixed type admitting the
data $(G, \tau)$, where $\tau=(g'|m_1, \cdots ,m_r)$. The number of connected components is equal
to the number of classes of generating vectors of type $\tau$
of $G^{\circ}$ modulo the action given by $Map_{g',[r]} \times
\Aut(G)$ where the first group acts via Hurwitz moves. In our case set $\mathbb{Z}_4=\langle x
\mid x^4=1 \rangle$, and since $G^{\circ}=\mathbb{Z}_2$, then the only generating vector of type $(2 \mid -)$ is given
by:
\[ a_1=x^2, \ \  b_1=1, \ \  a_2=1, \ \  b_2=1.
\]
This computation and the following proposition conclude the proof of
Theorem \ref{zerozero}.

%---------------------------------
\begin{prop}\label{irridcomp}  Each item in Table 2 provides exactly one irreducible subvariety of the moduli space of surfaces of general type.
\end{prop}
\begin{proof}
Recall from Theorem \ref{mainisot} that each item in Table 2 gives
rise to a surface $S$ of general type which is the minimal
desingularization of $T=(C \times F)/G$ and both $C/G$ and $F/G$
are elliptic curves.
\\ To see that each item in the Table 2 gives rise only to one topological type, we proceed analyzing case by case.  We have to prove that each pair of generating vectors is unique up to Hurwitz moves and simultaneous conjugation.
Hence for the first case there is nothing to prove.
\\ For the other cases the groups $G$ are all $(1 \mid m_1)-$generated, and denote by $a_1$,$b_1$ and $c_1$ the elements of a generating
vector ($|c_1|=m_1$). Recall  that the Hurwitz moves in this
case are generated by (see \cite[Corollary 1.11]{pol1}):
\[
\mathbf{1}: \left\{
\begin{array}{rl}
a_1 & \rightarrow a_1  \\
b_{1}  & \rightarrow  b_{1}a_{1} \\
c_1 & \rightarrow c_1 \\
\end{array}
\right.
\mbox{   }
\mathbf{2}: \left\{
\begin{array}{rl}
a_1 & \rightarrow a_1b^{-1}_1  \\
b_{1}  & \rightarrow  b_{1} \\
c_1 & \rightarrow c_1. \\
\end{array}
\right.
\]
Notice that in all the cases the groups $G$ have the property that
$[G,G]-\{id\}$ consists of a unique conjugacy class. Hence we can
fix $c_1 \in [G,G]$.
\\ In case $G=D_4$, let us fix a rotation $x$ and a reflection $y$. For what we said $c_1=x^2$.
Moreover we see that $a_1$ and $b_1$ cannot be both rotations, and up to Hurwitz moves we can
assume that are both reflections. The two reflections must also be in two different conjugacy
classes in order to generate $G$. Applying then simultaneous conjugation we see that the generating vector is unique.
\\ In case $G=Q_8$, $c_1=-1$ and since the vector must generate $G$, up to simultaneous conjugation
the pair $(a_1,b_1)$ is one of the following: $(i,j)$, $(j,i)$,
$(i,k)$, $(k,i)$, $(k,j)$, $(j,k)$. By Hurwitz moves all the pairs
are equivalent to the one: $a_1=i$ and $b_1=j$, hence the
generating vector is unique.
\\ In case $G=S_3$, $c_1=(123)$, and since the vector must generate $G$, $a_1$ and $b_1$ cannot be both $3-$cycles.
Moreover up to Hurwitz moves we can assume that both $a_1$ and $b_1$ are transpositions. Since all the transpositions in $S_3$ are conjugate,
we see that the generating vector is unique.
\\ In the last case $G=A_4$, $c_1=(12)(34)$. To generate $G$ we need a $3-$cycle, hence $a_1$ and $b_1$
cannot be both $2-2$-cycles. Up to Hurwitz moves we can suppose that both are $3-$cycles, which however might be in different
conjugacy classes. But again applying Hurwitz moves we can suppose that they are in the same conjugacy class, hence the generating vector is unique.
\\ In the end by Theorem 5.4 \cite{cat00} we have that isotrivial fibred surfaces with fixed topological
type form a union of irreducible subvarieties of the moduli space of surfaces of general type.
Here for each case we have only one irreducible variety, whose
dimension is $2$ for all the cases except the first where the
dimension of the variety is $4$. The calculation of the dimension is done in the same way as for surfaces isogenous to a product.
\end{proof}
%-------------------------------------------------------------------------------------%

Since our program is written in all generality, we can complete
the classification of the surfaces isogenous to a product with
$p_g=q=1$ given in \cite{canpol} adding the dimension of the
moduli spaces and the number of connected components.
\begin{theo}[\cite{canpol}]\label{polizzimoduli} Let $S= (C \times F)/G$ be a
surface with $p_g=q=1$ isogenous to a product of curve. Then $S$
is minimal of general type and the occurrence for $g(C)$, $g(F)$
and $G$ are precisely those in the Table 3 in the Appendix, where $dim$ is the
dimension of the moduli space, and $n$ is the number of connected
components.
\end{theo}
%--------------------------------------------------------------------------------------%
\section{Fundamental Groups}\label{fundgr}
To study the component of the moduli space relative to a surface $S$ it is sometimes useful to know the fundamental group of the surface in question.
In this section we will calculate the fundamental group of the found isotrivial surfaces.
In case of surface isogenous to a product of curves we have the following.
\begin{prop}[\cite{cat00}] Let $S:=(C_1 \times C_2)/G$ be isogenous to a product of curves. Then the fundamental group of $S$ sits in an exact sequence:
\[ 1 \longrightarrow \Pi_{g_1} \times \Pi_{g_2} \longrightarrow \pi_1(S) \longrightarrow G \longrightarrow 1,
\]
where $\Pi_{g_i}:=\pi_1(C_i)$, and this extension is determined by the associated maps $G \rightarrow Map_{g_i}$ to the respective Teichm\"uller modular groups.
\end{prop}
By \cite{ingefabfritzroby} there is a similar description of the
fundamental group of isotrivial surfaces too, which enables us to
describe the fundamental group of the surfaces of Table 2.
Following \cite{ingefabfritzroby} we have:
\begin{defin}Let $G$ be a finite group and let $\theta_i:\Gamma_i:=\Gamma(g'_i|m_1,\dots,m_{k_i}) \rightarrow G$ be two admissible epimorphisms. We define the fibre product $\mathbb{H}:=\mathbb{H}(G,\theta_1,\theta_2)$ by:
\[ \mathbb{H}:=\{(x,y) \in \Gamma_1 \times \Gamma_2 \mid \theta_1(x)=\theta_2(y)\}.
\]
In fact $\mathbb{H}$ is defined by the cartesian diagram:
\[
\begin{xy}
\xymatrix{%%
\mathbb{H} \ar[rr] \ar[dd] & & \Gamma_1 \ar[dd]^{\theta_1} \\ \\
 \Gamma_2 \ar[rr]^{\theta_2} & & G.
  }
\end{xy}
\]
\end{defin}
\begin{defin} Let $H$ be a group. Then its torsion subgroup $\mathbf{Tors}(H)$ is the (normal) subgroup generated by elements of finite order in H.
\end{defin}
\begin{prop}[Proposition 3.4 \cite{ingefabfritzroby}] Let $C_1,C_2$ be compact Riemann surfaces of respective
genera $g_i \geq 2$ and let $G$ be a finite group acting faithfully on each $C_i$ as a group of biholomorphic transformations. Let $T=(C_1 \times C_2)/G$, and denote by $S$ a minimal desingularization of $T$. Then the fundamental group $\pi_1(T) \cong \pi_1(S)$. Moreover let $\theta_i:\Gamma_i:=\Gamma(g'_i|m_1,\dots,m_{k_i}) \rightarrow G$ be two admissible epimorphisms. Then $\pi_1(T) \cong \mathbb{H}/\mathbf{Tors}(\mathbb{H})$.
\end{prop}
Moreover we have the following structure theorem in the hypothesis of the above theorem.
\begin{theo}[Theorem 0.3 \cite{ingefabfritzroby}] The fundamental group $\pi_1(S)$ has a normal subgroup $\mathcal{N}$ of finite index which is isomorphic to the product of surface groups, i.e., there are natural numbers $h_1,h_2 \geq 0$ such that $\mathcal{N} \cong \pi_1(\widehat{C}_1) \times \pi_1(\widehat{C}_2)$, where $\widehat{C}_i$ are smooth curves of genus $h_i$ respectively.
\end{theo}
We have then the following theorem.
\begin{theo} The fundamental group of the surfaces given by the first 4 items in Table 2 is $\mathbb{Z}^4$. While the fundamental group $P$ of the last surface fits into the exact sequence:
\[ 1 \longrightarrow \mathbb{Z}^2 \times \mathbb{Z}^2 \longrightarrow P \longrightarrow D_4 \curlyvee D_4 \longrightarrow 1,
\]
where $D_4 \curlyvee D_4$ is the central product of $D_4$ times $D_4$, which is an extraspecial group of order $32$.
\end{theo}
\begin{proof}
To compute a presentation of the fundamental group we use the program implemented in MAGMA given in \cite{ingefabfritzroby}, with few modifications for orbifold groups $\Gamma(g|m_1,\dots,m_k)$ with $g = 1$. The program gives us a presentation for all the groups. In the first $4$ cases of Table 2, one has $\mathbb{Z}^4$. The last case we have a presentation given by:
\\  \texttt{ Finitely presented group P on 4 generators Relations}
\[ [P_1, P_2] = Id(P) \]
\[ [P_3, P_2] = Id(P) \]
\[ [P_4, P_3] = Id(P) \]
\[   P^{-1}_2 * P^{-1}_1 * P^{-1}_4 * P_2 * P_1 * P_4 = Id(P) \]
\[ P_1 * P^{-1}_3 * P_4 * P^{-1}_1 * P^{-1}_4 * P_3 = Id(P)\]
\[ P^{-2}_1 * P_4 * P^2_1 * P^{-1}_4 = Id(P) \]
\[ P^{-1}_3 * P_1 * P_3 * P_4 * P^{-1}_1 * P^{-1}_4 = Id(P) \]
\[ P^{-1}_2 * P^2_4 * P_2 * P^{-2}_4 = Id(P) \]
Form this presentation one notices that the square of the generators all lie in the center of $P$.
The core $C$ of the subgroup $<P_1,P_2,P^2_3,P^2_4>$ is isomorphic to $\mathbb{Z}^4$ and the group $P/C$ is identify by MAGMA as the small group of order $32$ and MAGMA library-number $49$, which is the central product of $D_4$ times $D_4$, known as the extraspecial group $D_4 \curlyvee D_4$. After inspection one sees that $32$ is the minimal index.
\end{proof}
%-------------------------------------------------------------------------------%

\section[Appendix]{Appendix: The computation of the number of connected
components (S\"onke Rollenske)}
\newcommand{\margincom}[1]{\marginpar{\sffamily\small #1}}

In this appendix we want to explain the strategies used to compute the number
of connected components of the Moduli space of surfaces isogenous to a
product with fixed invariants. All computations have been performed using the
computer algebra system GAP \cite{GAP4}.

For the following discussion let us fix a group $G$ together with two types of
generating vectors $\tau_i:=(g_i\mid m_1, \dots, m_{k_i})$. We want to find
the number $\mu$ of components of the moduli space of surfaces isogenous to a
product which can be constructed from this data.
%\margincom{Is there a better expression, you used before?}

At first sight it seems a pretty simple thing to do since we have the
group-theoretical description from the result of Bauer and Catanese (see
Theorem \ref{Fabmain}) and could try the following na\"{\i}ve approach:
\begin{itemize}
 \item Calculate all possible generating vectors for $G$ of the given types.
\item Check the compatibility of these pairs, i.e., if the action of $G$ on
$C\times F$ is free, and form the set of all such possible pairs.
\item Calculate the orbits under the group action given in Theorem
\ref{Fabmain}.
\end{itemize}
Then every orbit corresponds to exactly one component of the moduli space.

But this simple approach  turns out not to be feasible for the following
reason: the orbits are simply to big to be calculated in a reasonable amount
of time (and memory).
To give a concrete example we can consider $G=GL_2(\mathbb{F}_3)$, which has
number $(48, 29)$ in GAP, $\tau_1=(2\mid -)$ and $\tau_2=(0\mid 2,3,8)$. Then
there are 59719680 pairs of generating vectors which give rise to a surface
isogenous to a product with $p_g=q=2$ but there is only 1 orbit, i.e., 1 component of the moduli space.

Therefore we need a refined strategy in order to avoid the calculation of the complete
orbits. To explain it we need some notation: let $\mathfrak V_i$ be the set
of generating vectors of type $\tau_i$ for $G$ and let $\mathfrak X\subset
\mathfrak V_1\times \mathfrak V_2$ be the set of all compatible pairs of
generating vectors.

Let $H$ be the subgroup of the group of permutations of  $\mathfrak V_1\times
\mathfrak V_2$ generated by the action of the mapping class groups on both
factors and the diagonal action of the automorphism group of $G$ as described
in Theorem \ref{Fabmain}. We denote by $H_i$ the restriction of the action of
$H$ to the  component $\mathfrak V_i$.

The action of $H$ restricts to $\mathfrak X$ and $\mu$ coincides with number
of orbits in $\mathfrak X$.

The following, mostly obvious, lemma allows us to greatly simplify the
calculations.
\begin{lem}\label{lemlem}
Let $M_i$ be the mapping class group acting on $\mathfrak V_i$ and let
$(\mathcal V_1, \mathcal V_2)$ and $(\mathcal W_1, \mathcal W_2)$ be two
pairs of generating vectors. Then
\begin{enumerate}
 \item If $\mathcal V_1$ and $\mathcal W_1$ lie in the same $M_1$-orbit and
$\mathcal V_2$ and $\mathcal W_2$ lie in the same $M_2$-orbit then $(\mathcal
V_1, \mathcal V_2)$ and $(\mathcal W_1, \mathcal W_2)$ lie in the same
$H$-orbit.
\item If $\mathcal V_1$ and $\mathcal W_1$ do not lie in the same $H_1$-orbit
then  $(\mathcal V_1, \mathcal V_2)$ and $(\mathcal W_1, \mathcal W_2)$ lie
in different $H$-orbits.
\end{enumerate}
\end{lem}

Thus our revised algorithm takes roughly the following form:
\begin{itemize}
 \item Calculate a set $\mathfrak R_1$ of representants  of the $H_1$-orbits
on $\mathfrak V_1$, the generating vectors of type $\tau_1$, and a set
$\mathfrak R_2$ of representants of $M_2$-orbits on $\mathfrak V_2$.
 \item After testing the pairs in $\mathfrak R_1\times \mathfrak R_2$ for
compatibility we obtain a set of pairs $\mathfrak R \subset \mathfrak X$.
Each orbit in $\mathfrak X$ contains at least 1 element in $\mathfrak R$ by
\ref{lemlem} (\textit{i}).
\item We already have some lower bound on the number of components: if
$(\mathcal V_1, \mathcal V_2), (\mathcal W_1, \mathcal W_2)\in \mathfrak R$
then, by \ref{lemlem} (\textit{ii}), they lie in different orbits if
$\mathcal V_1\neq \mathcal W_1$ or if $\mathcal V_2$ and $\mathcal W_2$  lie
in different $H_2$-orbits.
\item It remains to calculate the full orbit only in the following case: there
are  $(\mathcal V_1, \mathcal V_2), (\mathcal W_1, \mathcal W_2)\in \mathfrak
R$ such that  $\mathcal V_1=\mathcal W_1$ and $\mathcal V_2$ and $\mathcal
W_2$  lie in different $M_2$-orbits but in the same $H_2$-orbit.
\end{itemize}
The last step was only necessary in very few of the considered cases, so we
mostly could deduce the number of components without calculating a single
$H$-orbit in  $\mathfrak X$.

For the results of the calculations we refer to main text, namely Theorem
\ref{zerozero} and Theorem \ref{polizzimoduli}.

The reader interested in more technical details concerning the algorithm for the computation of the orbits or in a version of the program should contact us via email, for the script see also \cite{io}.

\subsection*{Acknowledgements} This research was partly carried out at
Imperial College London supported by a DFG Forschungsstipendium. I would like
to thank Matteo Penegini for his offer to join this project.
%\newpage
\begin{center}
Table 3
\end{center}
\begin{center}
\begin{tabular}{|c|c|c|c|c|c|c|}
  \hline
 $g(F)$ & $g(C)$ & $G$ & IdSmallGroup & $\mathbf{m}$ & $dim$ & $n$  \\\hline
 $3$ & $3$ &  $(\mathbb{Z}_2)^2$ $(^*)$& G(4,2) & $(2^2)$, $(2^6)$ & $5$ & $1$ \\\hline $3$ &
 $5$ &$(\mathbb{Z}_2)^3$ $(^*)$ & G(8,5) & $(2^2)$, $(2^5)$& $4$ & $1$ \\ \hline
 $3$ & $5$ & $\mathbb{Z}_2  \times \mathbb{Z}_4$ $(^*)$ & G(8,2) & $(2^2)$, $(2^2,4^2)$& $3$ & $2$   \\\hline
 $3$ & $9$ & $\mathbb{Z}_2 \times \mathbb{Z}_8$ $(^*)$  & G(16,5) & $(2^2)$, $(2,8^2)$ & $2$ & $1$ \\ \hline
 $3$ & $5$ &  $D_4$ & G(8,3) & $(2^2)$, $(2^2,4^2)$& $3$ & $1$  \\ \hline
 $3$ & $7$ &  $D_6$  & G(12,4)  & $(2^2)$, $(2^3,6)$& $3$ & $1$ \\ \hline
 $3$ & $9$ &  $\mathbb{Z}_2 \times D_4$ & G(16,11) & $(2^2)$, $(2^3,4)$ & $3$ & $1$ \\ \hline
 $3$ & $13$ &  $D_{2,12,5}$ & G(24,5) &$(2^2)$, $(2,4,12)$& $2$ & $1$ \\ \hline
 $3$ & $13$ & $\mathbb{Z}_{2} \times A_4$& G(24,13) &$(2^2)$, $(2,6^2)$& $2$ & $1$ \\ \hline
 $3$ & $13$ & $S_4$ & G(24,12)  & $(2^2)$, $(3,4^2)$& $2$ & $1$\\ \hline
 $3$ & $17$ &  $\mathbb{Z}_2 \ltimes (\mathbb{Z}_2 \times \mathbb{Z}_8)$ &  G(32,9) & $(2^2)$, $(2,4,8)$ &$2$ & $1$ \\ \hline
 $3$ & $25$ &  $\mathbb{Z}_2 \times S_4$ & G(48,48)& $(2^2)$, $(2,4,6)$ & $2$ & $1$ \\ \hline
 $4$ & $3$ & $S_3$ & G(6,1)& $(3)$,$(2^6)$ & $4$ & $1$  \\ \hline %WW
 $4$ & $5$ &  $D_6$ & G(12,4)& $(3)$, $(2^5)$ & $3$ & $1$  \\ \hline
 $4$ & $7$ & $\mathbb{Z}_3 \times S_3$ &  G(18,3) & $(3)$, $(2^2,3^2)$& $2$ & $2$  \\ \hline
 $4$ & $7$ & $\mathbb{Z}_3 \times S_3$ & G(18,3) & $(3)$, $(3,6^2)$ &$2$ & $1$  \\ \hline
 $4$ & $9$ & $S_4$ & G(24,12)& $(3)$, $(2^3,4)$ & $2$ & $1$  \\ \hline
 $4$ & $13$ & $S_3 \times S_3$ & G(36,10 )&$(3)$, $(2,6^2)$& $1$ & $1$  \\ \hline
 $4$ & $13$ & $\mathbb{Z}_6 \times S_3$ & G(36,12) &$(3)$, $(2,6^2)$& $1$ & $1$  \\ \hline
 $4$ & $13$ & $\mathbb{Z}_4\ltimes (\mathbb{Z}_3)^2$ & G(36,9) &$(3)$, $(2,4^2)$& $1$ & $2$  \\ \hline
 $4$ & $21$ & $A_5$ & G(60,5) & $(3)$,$(2,5^2)$& $1$ & $1$  \\ \hline
          \end{tabular}
\end{center}
\begin{center}
\begin{tabular}{|c|c|c|c|c|c|c|c|c|}
  \hline $g(F)$ & $g(C)$ & $G$ & IdSmallGroup & $\mathbf{m}$ & $dim$ & $n$  \\
  \hline
 $4$ & $25$ &  $\mathbb{Z}_3 \times S_4$ & G(72,42) & $(3)$, $(2,3,12)$& $1$ & $1$  \\ \hline
 $4$ & $41$ &  $S_5$ & G(120,34) & $(3)$, $(2,4,5)$& $1$ & $1$   \\ \hline
 $5$ & $3$ &  $D_4$ &  G(8,3) & $(2)$, $(2^6)$& $4$ & $1$  \\ \hline
 $5$ & $4$ &  $A_4$ & G(12,3) & $(2)$, $(3^4)$& $2$ & $2$ \\ \hline
 $5$ & $5$ &  $\mathbb{Z}_4 \ltimes (\mathbb{Z}_2)^2$ & G(16,3) & $(2)$, $(2^2,4^2)$ & $2$ & $3$  \\ \hline
 $5$ & $7$ &  $\mathbb{Z}_2 \times A_4$ & G(24,13) & $(2)$, $(2^2,3^2)$& $2$ & $2$  \\ \hline
 $5$ & $7$ &  $\mathbb{Z}_2 \times A_4$ & G(24,13) & $(2)$, $(3,6^2)$& $1$ & $1$  \\ \hline
 $5$ & $9$ &  $\mathbb{Z}_8 \ltimes (\mathbb{Z}_2)^2$ & G(32,5)& $(2)$, $(2,8^2)$ & $1$ & $1$  \\ \hline
 $5$ & $9$ &  $\mathbb{Z}_2 \ltimes D_{2,8,5}$ & G(32,7) & $(2)$, $(2,8^2)$& $1$ & $1$  \\ \hline
 $5$ & $9$ &  $\mathbb{Z}_4 \ltimes (\mathbb{Z}_4 \times \mathbb{Z}_2)$ & G(32,2)& $(2)$, $(4^3)$ & $1$ & $1$  \\ \hline
 $5$ & $9$ &  $\mathbb{Z}_4 \ltimes (\mathbb{Z}_2)^3$ & G(32,6) & $(2)$, $(4^3)$& $1$ & $1$  \\ \hline
 $5$ & $13$ & $(\mathbb{Z}_{2})^2 \times A_4$ & G(48,49) & $(2)$, $(2,6^2)$& $1$ & $1$  \\ \hline
 $5$ & $17$ & $\mathbb{Z}_4 \ltimes (\mathbb{Z}_2)^4$ & G(64,32) & $(2)$, $(2,4,8)$& $1$ & $1$ \\ \hline
 $5$ & $21$ & $\mathbb{Z}_5 \ltimes (\mathbb{Z}_{2})^4$ & G(80,49) & $(2)$, $(2,5^2)$& $1$ & $2$ \\ \hline
 $5$ & $5$ &  $D_{2,8,3} $ & G(16,8)  & $(2^2)$ & $2$ & $1$ \\ \hline
 $5$ & $5$ &  $D_{2,8,5} $ & G(16,6) & $(2^2)$ & $2$ & $2$  \\ \hline
 $5$ & $5$ &  $\mathbb{Z}_4 \ltimes (\mathbb{Z}_2)^2$ &  G(16,3)  & $(2^2)$ & $2$ & $1$ \\ \hline
\end{tabular}
\end{center}
In cases with $(^*)$ the dimension of the moduli space and the
number of connected components where already given in \cite{pol1}.
%-------------------------------------------------------------------------------%
%%%%%%%%%%%%%%%%%%%%%%%%%%%%%%%%%%%%%%%%%%%%%%%%%%%%%%%%%%%%%%%%%%%%%%%%%%%%%%%%%%%%
%
%
%%%%%%%%%%%%%%%%%%%%%%%%%%%%%%%%%%%%%%%%%%%%%%%%%%%%%%%%%%%%%%%%%%%%%%%%%%%%%%%%%%%%%

%%%%%%%%%%%%%%%%%%%%%%%%%%%%%%%%%%%%%%%%%%%%%%%%%%%%%%%%%%%5
 %%%%%%%%%%%%%%%%%%%%%%%%%%%%%%%%%%%%%%%%%%%%%%%%%%%%%%%%%%%5
\end{document}